\documentclass[american]{scrartcl}
\KOMAoptions{DIV=12, paper=a4, abstract=on}

\usepackage[utf8]{inputenc}
\usepackage[T1]{fontenc}
\usepackage{lmodern} 
\usepackage{babel, microtype, fixltx2e}

\usepackage{amsfonts, amsmath, amssymb, amsthm}

\usepackage{subcaption, pgfplots, tikz}
\makeatletter\g@addto@macro\@floatboxreset\centering\makeatother
\pgfplotsset{compat=newest,
   plot coordinates/math parser=false,
   scale only axis,
   xmajorgrids, xminorticks=false,
   ymajorgrids, yminorticks=false,
   every axis/.append style={font=\small, line width=0.75pt, mark = +,%
      scaled x ticks = false},
   legend style={at={(0.02,1.00)},anchor=north west,align=left,
      fill=none,draw=none,row sep=-0.4em},
}
\newcommand{\cWinslow}{green} \newcommand{\ltWinslow}{solid}
\newcommand{\cHuang}{blue}    \newcommand{\ltHuang}{dashed}
\newcommand{\cHR}{red}        \newcommand{\ltHR}{dotted}
\newcommand{\cN}{black}        \newcommand{\ltN}{solid}

\usepackage{hyperref} 
\hypersetup{%
   pdfauthor  ={W. Huang, L. Kamenski, R. D. Russell},
   pdftitle   ={A comparative numerical study of~meshing functionals
      for variational mesh adaptation},
   pdfsubject ={preprint}
}
\usepackage[capitalize]{cleveref} 
\crefformat{section}{#2\S#1#3}
\Crefformat{section}{#2Section~#1#3}
\crefmultiformat{section}{Sects.~#2#1#3}%
   { and~#2#1#3}{, #2#1#3}{ and~#2#1#3}
\crefformat{equation}{(#2#1#3)}
\Crefformat{equation}{Equation~(#2#1#3)}
\crefmultiformat{equation}{(#2#1#3)}%
   { and~(#2#1#3)}{, (#2#1#3)}{, and~(#2#1#3)}
\crefrangeformat{equation}{(#3#1#4) to~(#5#2#6)}

\newenvironment{keywords}%
   {\begin{trivlist}\item[]{\bfseries\sffamily Key words:} }
   {\end{trivlist}}
\newenvironment{AMS}%
   {\begin{trivlist}\item[]{\bfseries\sffamily AMS subject classifications:} }
   {\end{trivlist}}

\providecommand{\Abs}[1]{\left\lvert#1\right\rvert}
\providecommand{\V}[1]{\boldsymbol{#1}}
\providecommand{\dx}{\,d\V{x}}

\providecommand{\p}[2]{\frac{\partial{}#1}{\partial{}#2}}
\providecommand{\Th}{\mathcal{T}_h}

\providecommand{\M}{\mathbb{M}}

\providecommand{\J}{\mathbb{J}}
\providecommand{\R}{\mathbb{R}}
\providecommand{\JMJ}{\J \M^{-1}\J^T}

\DeclareMathOperator{\tr}{tr}

\theoremstyle{plain}

\theoremstyle{definition}

\theoremstyle{remark}
\newtheorem{example}{Example}[section]

\title{%
   A~comparative numerical study of~meshing functionals
      for~variational mesh adaptation%
   \thanks{%
   Supported in~part
         by~NSF (U.S.A.) through grant DMS-1115118
         and~NSERC (Canada) through grant A8781.%
   }%
}

\author{%
   Weizhang~Huang%
   \thanks{%
      Department of~Mathematics, University of~Kansas, Lawrence, KS~66045, U.S.A.
      (\href{mailto:whuang@ku.edu}{\nolinkurl{whuang@ku.edu}}).%
   }%
   \and
   Lennard~Kamenski%
   \thanks{%
      Weierstrass Institute for~Applied Analysis and~Stochastics, Berlin, Germany
      (\href{mailto:kamenski@wias-berlin.de}{\nolinkurl{kamenski@wias-berlin.de}}).%
   }%
   \and
   Robert~D.~Russell%
   \thanks{%
      Department of~Mathematics, Simon Fraser University, Burnaby, BC~V5A~1S6, Canada
      (\href{mailto:rdr@sfu.ca}{\nolinkurl{rdr@sfu.ca}}).%
   }%
}

\date{}

\begin{document}
\maketitle

\begin{abstract}
We present a comparative numerical study for three functionals used for variational mesh adaptation.
One of them is a generalization of Winslow's variable diffusion functional while the others are based on equidistribution and alignment.
These functionals are known to have nice theoretical properties and work well for most mesh adaptation problems either as a stand-alone variational method or combined within the moving mesh framework.
Their performance is investigated numerically in terms of equidistribution and alignment mesh quality measures.
Numerical results in 2D and 3D are presented.


\end{abstract}

\section{Introduction}

Variational mesh adaptation is an important type of mesh adaptation method and has received considerable attention from scientists and engineers;  e.g.,\ see the books~\cite{HR11,KS94,Lis99,TWM85} and references therein.
It also serves as the base of a number of commonly used adaptive moving mesh methods (e.g.,\ see~\cite{CH01,HRR94a,HR99,LTZ01}).
In the variational approach, an adaptive mesh is generated as the image of a reference mesh under a coordinate transformation and such a coordinate transformation is determined as a minimizer of a certain meshing functional.
A number of meshing functionals have been developed in the past (cf.\ the above mentioned books).
For example, Winslow~\cite{Win81} proposed an equipotential method based on variable diffusion.
Brackbill and Saltzman~\cite{BS82} developed a method by combining mesh concentration, smoothness, and orthogonality.
Dvinsky~\cite{Dvi91} used the energy of harmonic mappings as his meshing functional, while Knupp~\cite{Knu96} and Knupp and Robidoux~\cite{KR00} developed functionals based on the idea of conditioning the Jacobian matrix of the coordinate transformation. 
More recently, Huang~\cite{Hua01b} and Huang and Russell~\cite{HR11} proposed functionals based on the so-called equidistribution and alignment conditions.

With variational mesh adaptation, the mesh concentration is typically controlled through a scalar or a matrix-valued function, often referred to as the metric tensor or monitor function and defined based on some error estimates and/or physical considerations.
While most of the meshing functionals have been developed with physical or geometric intuitions and have various levels of success in the adaptive numerical solution of partial differential equations (PDEs) and other applications, there is only a limited understanding on how the metric tensor affects the behavior of the mesh.
An attempt to alleviate this lack of understanding was made by Cao et al.~\cite{CHR99} for a generalization of Winslow's variable diffusion functional.
They showed that a significant change in an eigenvalue of the metric tensor along the corresponding eigendirection (first increasing and then decreasing, or vice versa) will result in adaptation of coordinate lines along this direction, although this adaptation competes with far more complicated effects, including those from changes in eigenvectors and other eigenvalues and the effects of the shapes of the physical and computational domains and the mesh point distribution on the boundaries.
In~\cite{Hua01b,HR11}, two functionals have been developed based directly on the equidistribution and alignment conditions.
These two conditions provide a complete characterization of the mesh elements through the metric tensor~\cite{Hua01b}.
Minimizing the functionals leads to meshes which tend to satisfy the conditions in an integral sense.
Nevertheless, this characterization is only qualitative, and how closely the mesh satisfies the conditions depends on the boundary correspondence between the computational and physical domains and the mesh point distribution on the boundaries.
Thus, numerical studies, especially comparative ones, are useful, and often necessary, in understanding how the mesh adaptation for those meshing functionals is controlled precisely by the metric tensor.
There do exist a few comparative numerical studies for meshing functionals.
For example, a gallery of (adaptive and non-adaptive) meshes is given in~\cite{KS94} for a number of meshing functionals.
Some comparative meshes are given in~\cite{HR11} for the harmonic mapping functional~\cite{Dvi91} and the subsequent functional based on equidistribution and alignment~\cite{Hua01b}.

The main objective of this work is to present a comparative study for three of the most appealing meshing functionals, a generalization of Winslow's variable diffusion functional (cf.~\cref{win-2}) and two functionals based on equidistribution and alignment (cf.~\cref{huang-1,hr-1}).
They are selected because \cref{win-2,huang-1} have been known to work well for many problems (e.g.,\ see~\cite{BMRS02,BMR01,Hua01b,HR97b,HR99,LTZ01}) while \cref{hr-1} is similar to \cref{huang-1} and has some very nice theoretical properties (cf. \cref{SEC:functional:EA}).
Another motivation is to present some three dimensional numerical results for those functionals.
Critical for our study is to perform the substantial computations using the improved implementation of the variational methods introduced in~\cite{HK2014}.
In a sharp contrast to the situation in two dimensions, very little work has been done with variational mesh adaptation and adaptive moving mesh methods in three dimensions (e.g.,\ see~\cite{HR11,LTZ01}).
It is particularly interesting to see how the functionals perform in this case.

An outline of the paper is given as follows.
We describe the basic ideas of the variational mesh adaptation and its direct discretization (that is, first to discretize and then optimize) in~\cref{SEC:variational}.
In~\cref{SEC:functionals} we introduce the three functionals to be studied for the numerical comparison, a generalization of Winslow's variable diffusion functional and
two functionals based on equidistribution and alignment.
Numerical results and example adaptive meshes are given in~\cref{SEC:numerics}, followed by conclusions in~\cref{SEC:conclusions}.

\section{Variational mesh adaptation}
\label{SEC:variational}

In variational mesh adaptation, an adaptive mesh is generated as the image of a reference mesh under a coordinate transformation.
Denote the physical domain by $\Omega \subset \R^d$ ($d \ge 1$), and assume that we are given a computational domain $\Omega_c \subset \R^d$ and a quasi-uniform mesh $\hat{\Th}_c$ thereon (in this work we consider only simplicial meshes).
In many situations we can choose $\Omega_c$ to be the unit square/cube or simply $\Omega$.
Denote the coordinate transformation by $\V{x} = \V{x}(\V{\xi}) \colon \Omega_c \to \Omega$ and its inverse by $\V{\xi} = \V{\xi}(\V{x}) \colon \Omega \to \Omega_c$.
Such a coordinate transformation (more precisely, its inverse) is determined as a minimizer of a meshing functional.
Most of the existing meshing functionals can be cast in a general form as
\begin{equation}
   I[\V{\xi}] = \int_\Omega G(\J, \det(\J), \M, \V{x}) \dx,
   \label{fun-1}
\end{equation}
where $G$ is a smooth function, $\J$ is the Jacobian matrix of $\V{\xi} = \V{\xi}(\V{x})$, $\det(\J)$ denotes the determinant of $\J$, and $\M = \M(\V{x})$ is the metric tensor supplied by the user to control mesh concentration.
We assume that $\M = \M(\V{x})$ is a symmetric and uniformly positive definite $d$-by-$d$ matrix-valued function on $\Omega$.
Notice that \cref{fun-1} is formulated in terms of the inverse coordinate transformation.
One reason for this is that this form is less likely to produce singular meshes~\cite{Dvi91}.
Another reason is that $\M$ is a function of $\V{x}$ and thus finding the functional derivative of $I[\V{\xi}]$ will not directly involve the derivatives of $\M$.
This is convenient since in practice $\M$ is known only at the vertices of a mesh and its derivatives are not cheap to find.
The main disadvantage of the formulation in this form is that $\V{\xi} = \V{\xi}(\V{x})$ (or its numerical approximation) does not give the physical mesh directly.
This is remedied either by interchanging the roles of the independent and dependent variables in the Euler-Lagrange equation of $I[\V{\xi}]$ (e.g.,\ see~\cite{HR11}) or, in a recently developed implementation (see below), computing the new physical mesh from a computational one using linear interpolation.

A minimizer of \cref{fun-1} can be found numerically in the MMPDE (moving mesh PDE) framework.
A conventional implementation~\cite{HR11} is to find the functional derivative of \cref{fun-1} and then define the MMPDE as the gradient flow equation of the functional. Having been transformed by interchanging the roles of the dependent and independent variables, the MMPDE can be discretized on $\hat{\Th}_c$ and a system of equations for the nodal velocities is obtained.
Finally, the new mesh is obtained by integrating the mesh equation over a time step.

A much simpler implementation was proposed recently by Huang and Kamenski~\cite{HK2014}.
Instead of utilizing the MMPDE directly, the new implementation first discretizes the functional on the current mesh $\Th$ and then, following the idea of the MMPDE, defines the mesh equation as the gradient equation of the discretized functional (with respect to the computational coordinates of the vertices).
To be specific, denote by $\V{x}_i$, $\hat{\V{\xi}}_i$, and $\V{\xi}_i$, $i= 1, \dotsc, N_v$ the coordinates of the vertices of the current physical mesh ($\Th$), the reference mesh ($\hat{\Th}_c$), and the computational mesh (${\Th}_c$), respectively.
We assume that these meshes have the same numbers of the elements and vertices and the same connectivity.
For any element $K \in \Th$ (with vertices $\V{x}_i^K$, $i = 0, \dotsc, d$), the corresponding element in ${\Th}_c$ is denoted by $K_c$ (with vertices $\V{\xi}_i^K$, $i=0, \dotsc, d$).
The edge matrices for $K$ and $K_c$ are defined as
\[
   E_K = [\V{x}_1^K-\V{x}_0^K, \dotsc, \V{x}_d^K-\V{x}_0^K]
   ,
   \quad
   E_{K_c} = [\V{\xi}_1^K-\V{\xi}_0^K, \dotsc, \V{\xi}_d^K-\V{\xi}_0^K]
   .
\]
Let $\omega_i$ be the element patch associated with vertex $\V{x}_i$ (i.e.,\ the collection of the elements containing $\V{x}_i$ as a vertex).
Then, the equation for the nodal velocities reads as
\begin{equation}
   \begin{cases}
   \frac{d \V{\xi}_i}{d t} 
      = \frac{P_i}{\tau} \sum\limits_{K \in \omega_i} |K| \V{v}_{i_K}^K,
      & \quad i = 1, \dotsc, N_v, \quad t_n < t \le t_{n+1},
   \\
   \V{\xi}_i(t_n) = \hat{\V{\xi}}_i,
   & \quad i = 1, \dotsc, N_v
   ,
   \end{cases}
   \label{mmpde-1}
\end{equation}
where  $|K|$ is the volume of $K$, $\V{v}_{i_K}^K$ is the local mesh velocity associated with vertex $\V{x}_i$ in $K$, $i_K$ denotes the local index of $\V{x}_i$ in $K$, $\tau > 0$ is a constant parameter used to adjust the time scale of mesh movement, and $P= (P_1, \dotsc, P_{N_v})$ is a positive function used to make the mesh equation
to have desired invariance properties.
Although the parameter $\tau$ can be absorbed in $P$, using two parameters has the advantage that the role of each parameter is clear: $\tau$ for the time scale of mesh movement while $P$ for the invariance properties.
Ideally, $\tau$ should be chosen such that the mesh movement has the same scale as the physical equation.
Unfortunately, there is no theoretical analysis for this yet and trial and error is still the most practical way to choose $\tau$.
Numerical experience shows that a value in the range $[0.01, 0.1]$ seems to work well for most problems~\cite{HR11}.
In our computation, we choose $P$ such that the equation is invariant under the scaling transformation $\M \to c\, \M$ for all non-zero constants $c$ (cf.~\cref{SEC:functionals}): in variational mesh adaptation it is the relative distribution of $\M$ over
the physical domain (instead of its absolute distribution) that determines the variation of the mesh density and therefore it is essential for the moving mesh equation to be invariant under the scaling transformation of $\M$.

The local velocities are given by
\begin{equation}
   \begin{bmatrix} (\V{v}_1^K)^T\\ \vdots \\  (\V{v}_d^K)^T \end{bmatrix} 
      = - E_K^{-1} \frac{\partial G}{\partial \J}
         - \frac{\partial G}{\partial \det(\J)}
            \frac{\det(E_{K_c})}{\det(E_K)} E_{K_c}^{-1}
   ,
   \qquad
   \V{v}_0^K =  - \sum_{j=1}^d \V{v}_j^K
   ,
   \label{mmpde-2}
\end{equation}
where the derivatives of $G$ with respect to $\J$ and $\det(\J)$ (see \cref{win-3,huang-2,hr-2} below) are evaluated as
\begin{align*}
   \frac{\partial G}{\partial \J} 
      & = \frac{\partial G}{\partial \J}
         \left(E_{K_c}  E_K^{-1}, 
            \frac{\det(E_{K_c})}{\det (E_K)}, \M(\V{x}_K), \V{x}_K\right),
   \\
   \frac{\partial G}{\partial \det(\J)} 
      & = \frac{\partial G}{\partial \det(\J)}
         \left(E_{K_c}  E_K^{-1}, 
            \frac{\det(E_{K_c})}{\det (E_K)}, \M(\V{x}_K), \V{x}_K\right)
   .
\end{align*}
The second equation in \cref{mmpde-2} is an inherent property
resulting directly from the differentiation of the discretized meshing functional; it states that the centroid of $K$ stays fixed if only the contribution from $K$ is taken into account. 

The above mesh equation should be modified properly for boundary vertices.
For example, if $\V{\xi}_i$ is a fixed boundary vertex, we replace the corresponding equation by
\[
   \p{\V{\xi}_i}{t} = 0.
\]
When $\V{\xi}_i$ is allowed to move on a boundary curve/surface represented by
\[
   \phi(\V{\xi}) = 0,
\]
then the mesh velocity $\p{\V{\xi}_i}{t}$ needs to be modified such that its normal component along the curve or surface is zero, i.e.,
\[
   \nabla \phi (\V{\xi}_i) \cdot \p{\V{\xi}_i}{t} = 0.
\]

Mesh equation \cref{mmpde-1} defines the movement of the nodes of the computational mesh ${\Th}_c$ starting from the reference mesh $\hat{\Th}_c$ at $t_n$.
The equation can be integrated in time to obtain the computational mesh at $t_{n+1}$.
For notational simplicity, we denote the computational mesh at $t_{n+1}$ by ${\Th}_c$ as well.
Notice that the physical mesh $\Th$ is fixed during the time integration from $t_n$ to $t_{n+1}$.
Meshes ${\Th}_c$ and $\Th$ define a correspondence 
\[
   \Th = \Psi({\Th}_c) .
\]
The new physical mesh is computed by means of linear interpolation as
\[
   \tilde{\Th} = \Psi(\hat{\Th}_c)
   ,
\]
where $\hat{\Th}_c$ is the reference mesh on $\Omega_c$.
The computational mesh
plays the role of an intermediate variable.

Recall that the mesh concentration in variational mesh adaptation is controlled through the metric tensor $\M = \M(\V{x})$.
Such a metric tensor can be defined based on physical or geometric considerations or some error estimates.
For example, for the $L^2$ norm of the error of piecewise linear interpolation on simplicial meshes, the optimal metric tensor~\cite{Hua05b,HS03} (also see~\cite[(5.192)]{HR11}) is
\begin{equation}
   \M = {\det \left(\alpha I +  \Abs{H(u)} \right)}^{- \frac{1}{d+4}} 
   \left [ \alpha I  +  \Abs{H(u)} \right ]
   ,
   \label{M-1}
\end{equation}
where $H(u)$ is the Hessian of function $u$, $|H(u)|$ is the eigen-decomposition of $H(u)$ with the eigenvalues being replaced by their absolute values, and the regularization parameter $\alpha > 0$ is chosen such that
\[
   \int_\Omega {\det(\M)}^{\frac{1}{2}} \dx 
   \equiv \int_\Omega {\det\left(\alpha I + |H(u)|\right)}^{\frac{2}{d+4}} \dx
   = 2 \int_\Omega {\det\left( |H(u)|\right)}^{\frac{2}{d+4}} \dx .
\]

In practical computation, $u$ is typically unknown, and only the approximations to its values at the vertices are available.
For this reason (and even in situations where an analytical expression for $u$ is available), the Hessian in \cref{M-1} is replaced by an approximation obtained by a Hessian recovery technique from the nodal values of $u$ or the approximations of the nodal values of $u$. 
A number of such techniques are known to produce nonconvergent recovered Hessians from a linear finite element approximation  (e.g.,\ see Kamenski~\cite{Kam09}).
Nevertheless, it is shown by Kamenski and Huang~\cite{KaHu2013} that a linear finite element solution of an elliptic BVP converges at a second order rate as the mesh is refined if the recovered Hessian used to generate the adaptive mesh satisfies a closeness assumption.
Numerical experiment shows that this closeness assumption is satisfied by the approximate Hessian obtained with commonly used Hessian recovery methods.
We use a Hessian recovery method based on a least squares fit: a quadratic polynomial is constructed locally for each vertex via least squares fitting to neighboring nodal function values and an approximate Hessian at the vertex is then obtained by differentiating the polynomial.

\section{Meshing functionals}
\label{SEC:functionals}

Here we introduce the three meshing functionals used in the numerical study.
A generalization of Winslow's variable diffusion functional and the two functionals based on equidistribution and alignment are selected because they are reasonably simple, have nice theoretical properties, and are known to work well for many problems.

\subsection{Winslow's functional based on variable diffusion}
\label{SEC:functional:win}

The first functional is the variable diffusion proposed by Winslow~\cite{Win81}. 
It uses the system of elliptic PDEs
\[
   - \nabla \cdot \left( w \nabla \xi_i \right) = 0, \quad i = 1, \dotsc, d
   ,
\]
for generating adaptive meshes, where $w = w(\V{x}) > 0$ is the weight function.
This system mimics a (steady-state) diffusion process with a heterogeneous diffusion coefficient $w(\V{x})$.
It is the Euler-Lagrange equation of the functional
\begin{equation}
   I[\V{\xi}] = \frac{1}{2} 
      \int_\Omega \sum\limits_{i=1}^d w(\V{x}) |\nabla \xi_i |^2 \dx
   = \frac{1}{2} \int_\Omega w(\V{x}) \tr (\J \J^T) \dx
   ,
   \label{win-1}
\end{equation}
where $\tr(\cdot)$ is the trace of a matrix.
A generalization of this functional to allow a diffusion tensor is
\begin{equation}
   I[\V{\xi}] = \frac{1}{2} \int_\Omega \tr (\J \M^{-1} \J^T) \dx
   .
   \label{win-2}
\end{equation}
This functional has been used by a number of researchers, e.g.,\ see Huang and Russell~\cite{HR97b,HR99}, Li et al.~\cite{LTZ01}, and Beckett et al.~\cite{BMR01}.
It is coercive and convex~\cite[Example 6.2.1]{HR11}.
Thus, under  a suitable boundary condition (such as the Dirichlet boundary condition with $\partial \Omega_c$ being mapped onto $\partial \Omega$), the functional \cref{win-2} has a unique minimizer.

For this functional, the derivatives of $G$ with respect to $\J$ and $\det(\J)$ needed in \cref{mmpde-2} are
\begin{equation}
   \begin{cases}
   G = \frac{1}{2} \tr (\J \M^{-1} \J^T),
   \\
   \frac{\partial G}{\partial \J} = \M^{-1} \J^T,
   \\
   \frac{\partial G}{\partial \det(\J)} = 0
   .
   \end{cases}
\label{win-3}
\end{equation}
The interested reader is referred to~\cite{HK2014} for the derivation.

The balancing function in \cref{mmpde-1} is chosen as $P = \det(\M)^{\frac{1}{d}}$.
With this choice, \cref{mmpde-1} is invariant under the scaling transformation $\M \to c\, \M$.

\subsection{Functionals based on equidistribution and alignment}
\label{SEC:functional:EA}

The other functionals are based on the equidistribution and alignment conditions.
These conditions provide a full mathematical characterization of a non-uniform mesh.
Indeed, any non-uniform mesh can be viewed as a uniform one in the metric specified by a tensor.
Moreover, a mesh is uniform in the metric specified by the metric tensor $\M = \M(\V{x})$ if and only if it satisfies the equidistribution and alignment conditions associated with $\M$~\cite{Hua06,HR11}.
In the continuous form, they are
\begin{align}
   \text{equidistribution:}\quad 
      & \det(\J)^{-1} \det(\M)^{\frac{1}{2}}
         = \frac{\sigma}{|\Omega_c|}, \quad \forall \V{x} \in \Omega
   \label{eq-1}
   \\
   \text{alignment:}\quad
   & \frac{1}{d} \tr (\J \M^{-1} \J^T) 
      = \det (\J \M^{-1} \J^T)^{\frac{1}{d}}, \quad \forall \V{x} \in \Omega,
   \label{ali-1}
\end{align}
where
\begin{equation}
   \sigma = \int_{\Omega}  \det(\M)^{\frac{1}{2}} \dx.
   \label{sigma-1}
\end{equation}
These conditions require the mesh elements to have the same size (equidistribution) and be equilateral (alignment) in the metric $\M$, respectively.
The alignment condition also implies that the elements are aligned with $\M$ in the sense that the principal directions of the circumscribed ellipsoid of each element coincide with the eigen-directions of $\M$ while the lengths of the principal axes of the ellipsoid are reciprocally  proportional to the square roots of the eigenvalues of $\M$. 

The first functional based on equidistribution and alignment, proposed in~\cite{Hua01b}, is
\begin{equation}
   I[\V{\xi}]  
   = \theta \int_\Omega \sqrt{\det(\M)} {\left(\tr(\JMJ) \right)}^{\frac{dp}{2}} \dx
   + (1 - 2 \theta) d^{\frac{dp}{2}} \int_\Omega \sqrt{\det(\M)}
            {\left( \frac{\det(\J)}{\sqrt{\det(\M)}}\right)}^{p} \dx,
   \label{huang-1}
\end{equation}
where $\theta \in (0,1)$ and $p > 0$ are dimensionless parameters.
Loosely speaking, the first and second terms correspond to the equidistribution and alignment conditions, respectively.
The terms are dimensionally homogeneous and the balance between them is controlled by the dimensionless parameter $\theta$.
For $0 < \theta \le \frac{1}{2}$, $d p \ge 2$, and $p \ge 1$, the functional is coercive and polyconvex and has a minimizer~\cite[Example 6.2.2]{HR11}.
Moreover, for $\theta = \frac{1}{2}$ and $d p = 2$ it reduces to
\[
   I[\V{\xi}]  = \frac{1}{2} \int_\Omega \sqrt{\det(\M)} \tr(\JMJ) \dx,
\]
which is exactly the energy functional of a harmonic mapping from $\Omega$ to $\Omega_c$ (cf.~\cite{Dvi91}).

For the functional \cref{huang-1}, we have
\begin{equation}
\begin{cases}
   G  = \theta \sqrt{\det(\M)} {\left(\tr(\JMJ) \right)}^{\frac{dp}{2}}
       + (1 - 2 \theta) d^{\frac{dp}{2}} \sqrt{\det(\M)}
         {\left( \frac{\det(\J)}{\sqrt{\det(\M)}}\right)}^{p},
   \\
   \p{G}{\J}  =  d p \theta \sqrt{\det(\M)}  
      {\left( \tr(\JMJ )\right)}^{\frac{d p}{2}-1} \M^{-1} \J^T,
   \\
   \p{G}{\det(\J)}  = p (1-2\theta) d^{\frac{d p}{2}} 
      \det{(\M)}^{\frac{1-p}{2}} \det{(\J)}^{p-1}
   .
   \end{cases}
   \label{huang-2}
\end{equation}
In the computation, we use $(p, \theta) = (2,\frac{1}{3})$.
$p = 2$ is the smallest integer to satisfy $d p \ge 2$ for $d = 1$, 2, and 3. 
The choice of $\theta = 1/3$ is in the range $(0, 1/2]$ for the functional to be polyconvex while giving a bigger weight to the equdistribution condition.
This set of the values works well for all tested problems.
The balancing function in \cref{mmpde-1} is chosen to be $P = \det(\M)^{\frac{p-1}{2}}$, so that \cref{mmpde-1} is invariant under the scaling transformation $\M \to c\ \M$.

The second functional based on equidistribution and alignment is
\begin{multline}
I[\V{\xi}] 
   = \theta_1 \int_\Omega \sqrt{\det(\M)} {\left(\tr(\JMJ) \right)}^{\frac{dp}{2}} \dx
   \\
   + \theta_2 d^{\frac{d p (d-2)}{2 (d-1)}} 
      \int_\Omega \det(\M)^{\frac{1}{2} (1 - \frac{d p}{d-1})}
      \det(\J)^{\frac{d p}{d-1}} 
      {\left(\tr(\J^{-T}\M \J^{-1}) \right)}^{\frac{d p}{2(d-1)}} \dx
   \\
   + \left(\theta_3 - \theta_1 - \theta_2 \right) d^{\frac{dp}{2}} 
      \int_\Omega \sqrt{\det(\M)}
      {\left( \frac{\det(\J)}{\sqrt{\det(\M)}}\right)}^{p} \dx       
   \\
   + \frac{\theta_4}{\sigma^{p+\nu}} \int_\Omega   \sqrt{\det(\M)}
      {\left( \frac{\det(\J)}{\sqrt{\det(\M)}}\right)}^{-\nu} \dx ,
   \label{hr-1}
\end{multline}
where $p > 1$, $\nu > 0$, and $\theta_i > 0$ ($i=1, \dotsc, 4$) are parameters.
The first three terms in \cref{hr-1} are dimensionally homogeneous in $\M$ and $\J$ while the last term has the same dimension in $\M$ as the other terms.
This functional was proposed in~\cite[(6.120)]{HR11} to avoid singularity of the coordinate transformation.
Indeed, if $\theta_3 - \theta_1 - \theta_2 > 0$, then the functional is coercive and polyconvex and has a minimizer satisfying $\det(\J) > 0$ in $\Omega$~\cite[Example~6.2.3]{HR11}.

In the computation, we choose $\theta_1 = \theta_2 = \frac{1}{3}$, $\theta_3 = 1$, $\theta_4 = 0.1$, $p = 2$, $\nu = 1$, and the balancing function $P = \det(\M)^{\frac{p-1}{2}}$.
These choices are based on the functional \cref{huang-1} and the desire to keep the fourth term relatively small.

For this functional, we have
\begin{equation}
\begin{cases}
   G 
   = \theta_1 \sqrt{\det(\M)} {\left(\tr(\JMJ) \right)}^{\frac{dp}{2}}
   \\ \qquad
   + \theta_2 d^{\frac{dp(d-2)}{2(d-1)}} \det(\M)^{\frac{1}{2} (1 - \frac{d p}{d-1})}
	   \det(\J)^{\frac{d p}{d-1}} 
      {\left(\tr(\J^{-T}\M \J^{-1}) \right)}^{\frac{d p}{2(d-1)}}
   \\ \quad \qquad
   + (\theta_3 - \theta_1 - \theta_2) d^{\frac{dp}{2}} 
      \left(\sqrt{\det(\M)}\right)^{1-p} \det(\J)^{p}
   \\ \qquad \qquad
   + \frac{\theta_4}{\sigma^{p+\nu}} 
      \left(\sqrt{\det(\M)}\right)^{1+\nu}\det(\J)^{-\nu} ,
   \\
   \frac{\partial G}{\partial \J}
   = \theta_1 d p \sqrt{\det(\M)} {\left(\tr(\JMJ) \right)}^{\frac{dp}{2}-1}\M^{-1}\J^T
   \\ \qquad
   - \frac{\theta_2 d p}{d-1} d^{\frac{d p (d-2)}{2 (d-1)}} \det(\M)^{\frac{1}{2}
      (1 - \frac{d p}{d-1})}
	   \det(\J)^{\frac{d p}{d-1}}
      {\left(\tr(\J^{-T}\M \J^{-1}) \right)}^{\frac{d p}{2(d-1)}-1}
      \J^{-1} \J^{-T} \M \J^{-1} ,
   \\
   \frac{\partial G}{\partial \det(\J)} 
   =  \frac{\theta_2 d p}{d-1} d^{\frac{d p (d-2)}{2 (d-1)}}
      \det(\M)^{\frac{1}{2} (1 - \frac{d p}{d-1})}
		\det(\J)^{\frac{d p}{d-1}-1}
      {\left(\tr(\J^{-T}\M \J^{-1}) \right)}^{\frac{d p}{2(d-1)}}
   \\ \qquad
   + (\theta_3 - \theta_1 - \theta_2) p d^{\frac{dp}{2}} 
      \left(\sqrt{\det(\M)}\right)^{1-p} \det(\J)^{p-1}
   \\ \quad \qquad
   - \frac{\theta_4\nu }{\sigma^{p+\nu}} 
      \left(\sqrt{\det(\M)}\right)^{1+\nu}\det(\J)^{-\nu-1} .
\end{cases}
\label{hr-2}
\end{equation}

\section{Numerical experiments}
\label{SEC:numerics}

In the following we consider a number of examples in two and three dimensions.
For a given function we consider $\M$ defined in \cref{M-1} which is optimal for minimizing the $L^2$ norm of the linear interpolation error of this function and compare meshes obtained from using the functionals of Winslow \cref{win-2} (W), Huang \cref{huang-1} (H), and Huang and Russell \cref{hr-1} (HR).

To assess the quality of the generated meshes, we compare the $L^2$ norm of the linear interpolation error and the equidistribution and alignment mesh quality measures, which describe how far the mesh is from being uniform in the metric defined by $\M$.
The element-wise quality measures are based on \cref{eq-1} and \cref{ali-1}  and defined as 
\begin{equation}
   Q_{eq,K} =
      \frac{ \det(\J_K)^{-1} \det(\M_K)^{\frac{1}{2}} }
           { \sigma / |\Omega_c|} ,
   \qquad  
   Q_{ali,K} = 
      \frac{ \tr (\J_K \M_K^{-1} \J_K^T)}
      {d \det(\J_K \M_K^{-1} \J_K^T)^{\frac{1}{d}}}
   ,
   \label{eq:mesh:quality}
\end{equation}
while for the overall mesh quality measures we take their root-mean-squared values,
\begin{equation}
   Q_{eq}  = \sqrt{ \frac{1}{N} \sum_{K\in\Th} Q_{eq, K}^2 } ,
   \qquad
   Q_{ali} = \sqrt{ \frac{1}{N} \sum_{K\in\Th} Q_{ali, K}^2 }
   .
   \label{eq:mesh:global:quality}
\end{equation}
If the mesh is uniform with respect to $\M$, then $Q_{eq} = Q_{ali} = 1$; if the mesh is far from being uniform with respect to $\M$, then $Q_{eq}$ and $Q_{ali}$ will become large.
In other words, these quality measures describe how well the volume (measured by $Q_{eq}$) and the shape and orientation (measured by $Q_{ali}$) of mesh elements correspond to the desired size and shape prescribed by $\M$ (see~\cite{Hua05c} for more details on the mesh quality measures).

\subsection{Two dimensions}
\begin{figure}[p]
   \begin{subfigure}[t]{1.0\linewidth}
      \includegraphics[clip, width=0.30\linewidth]{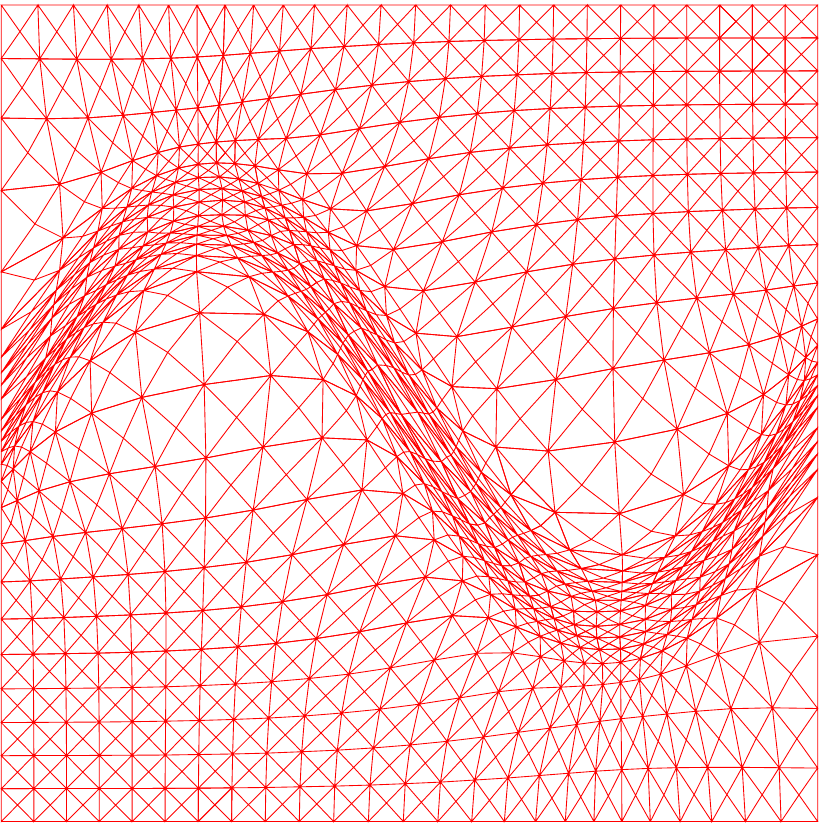}%
      \hfill%
      \includegraphics[clip, width=0.30\linewidth]{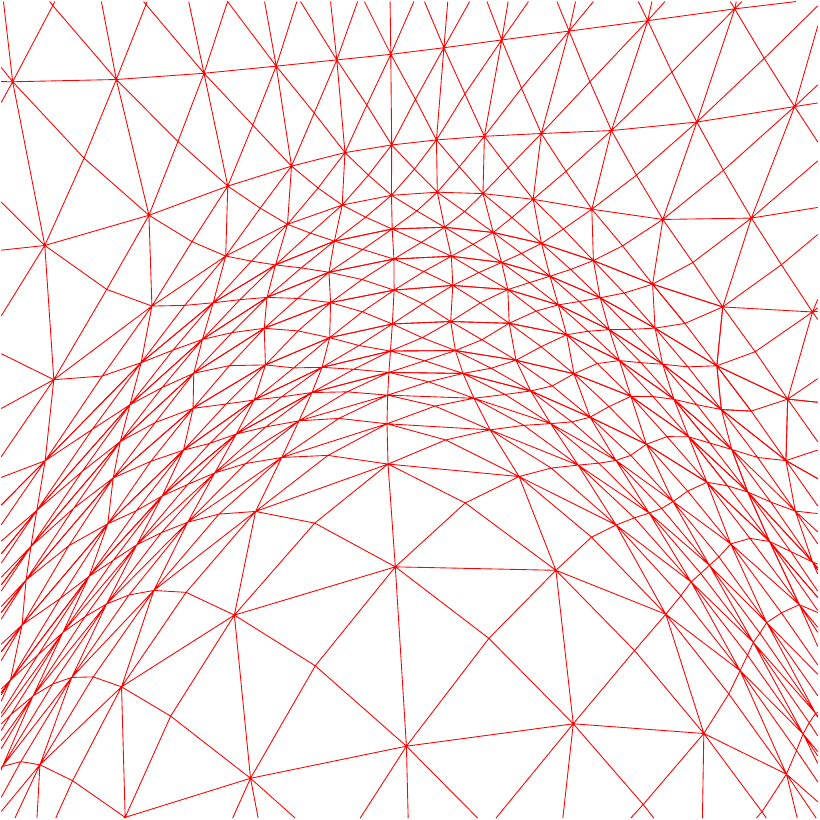}%
      \hfill%
      \includegraphics[clip, width=0.30\linewidth]{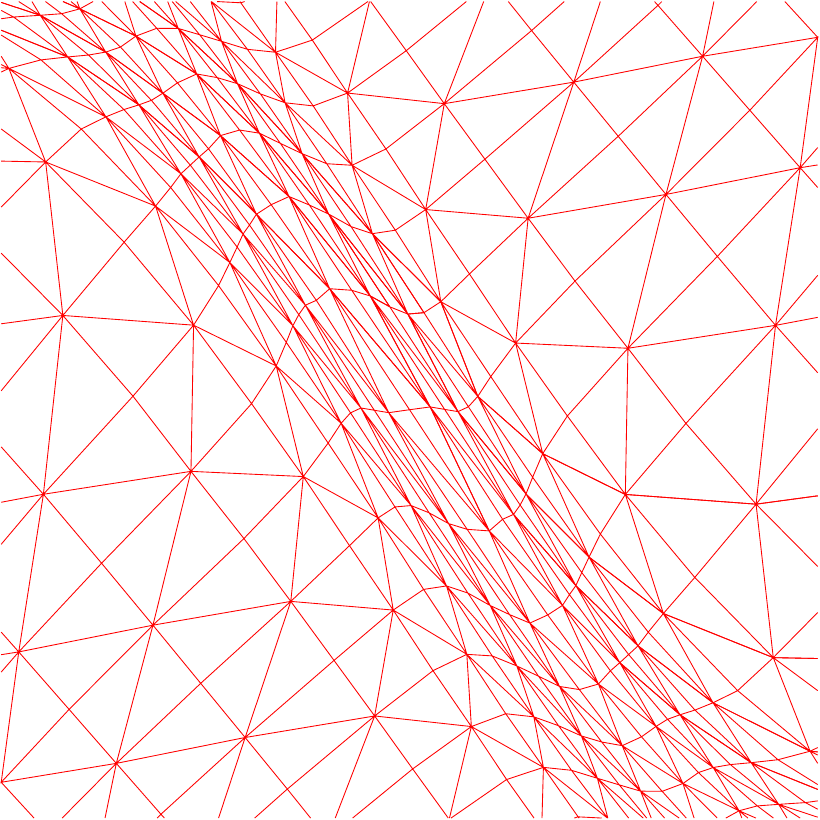}%
      \caption{Winslow's \cref{win-2}}\label{fig:ex:1:winslow}
   \end{subfigure}%
   \\[2.0ex]
   \begin{subfigure}[t]{1.0\linewidth}
      \includegraphics[clip, width=0.30\linewidth]{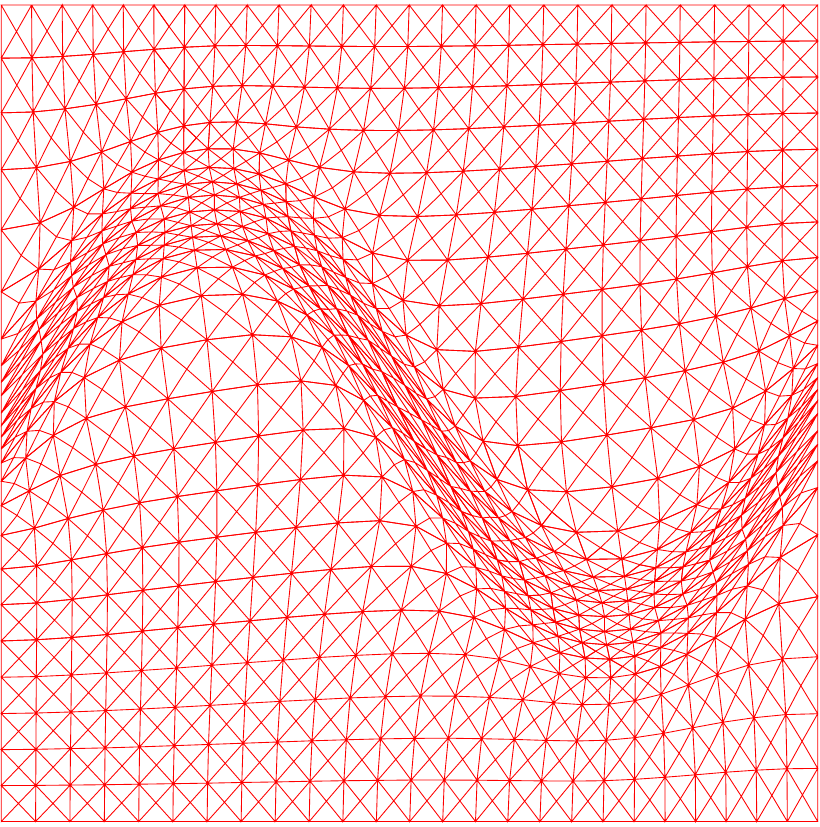}%
      \hfill%
      \includegraphics[clip, width=0.30\linewidth]{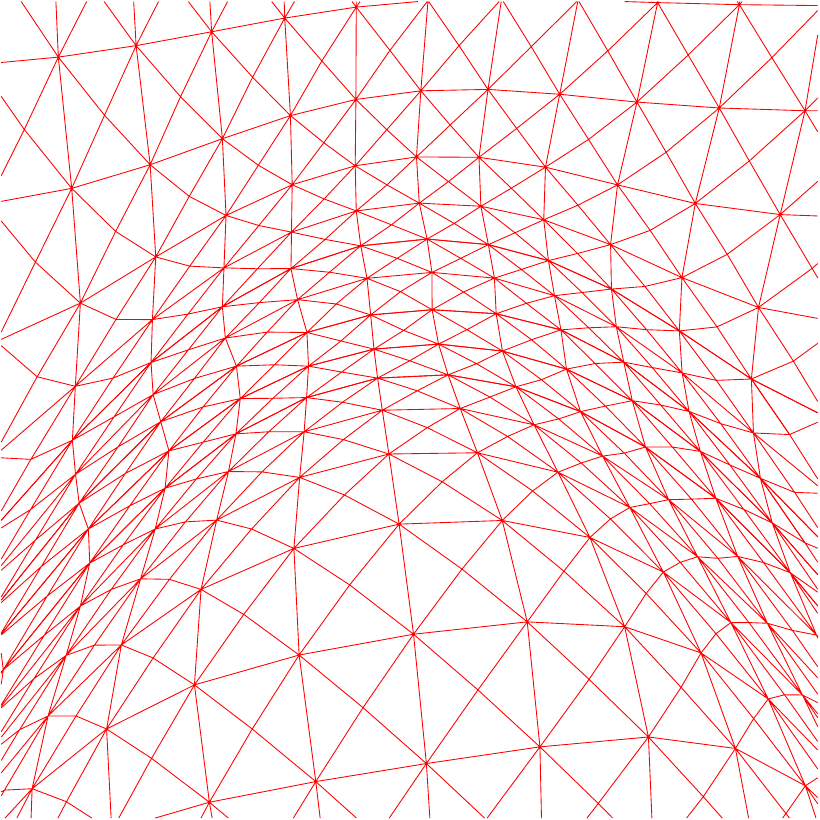}%
      \hfill%
      \includegraphics[clip, width=0.30\linewidth]{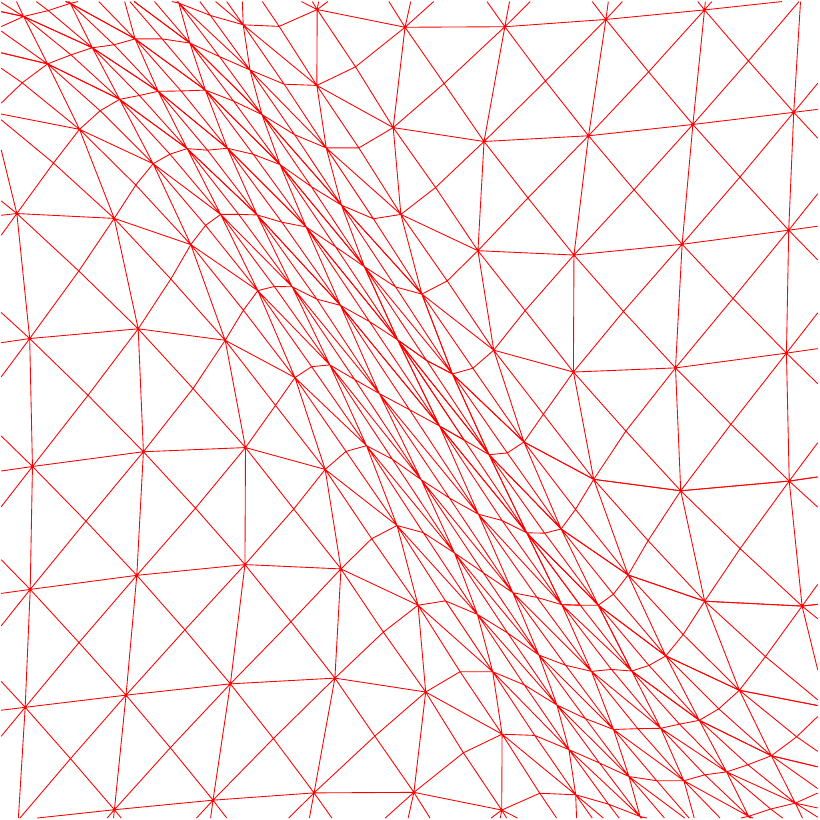}%
      \caption{Huang's \cref{huang-1}}\label{fig:ex:1:huang}
   \end{subfigure}%
   \\[2.0ex]
   \begin{subfigure}[t]{1.0\linewidth}
      \includegraphics[clip, width=0.30\linewidth]{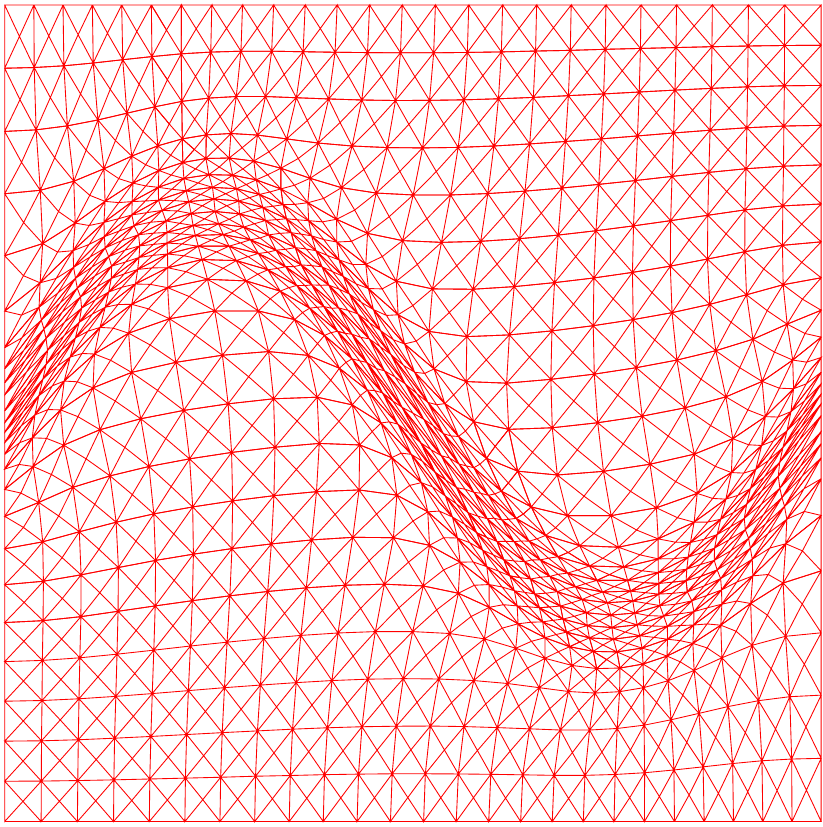}%
      \hfill%
      \includegraphics[clip, width=0.30\linewidth]{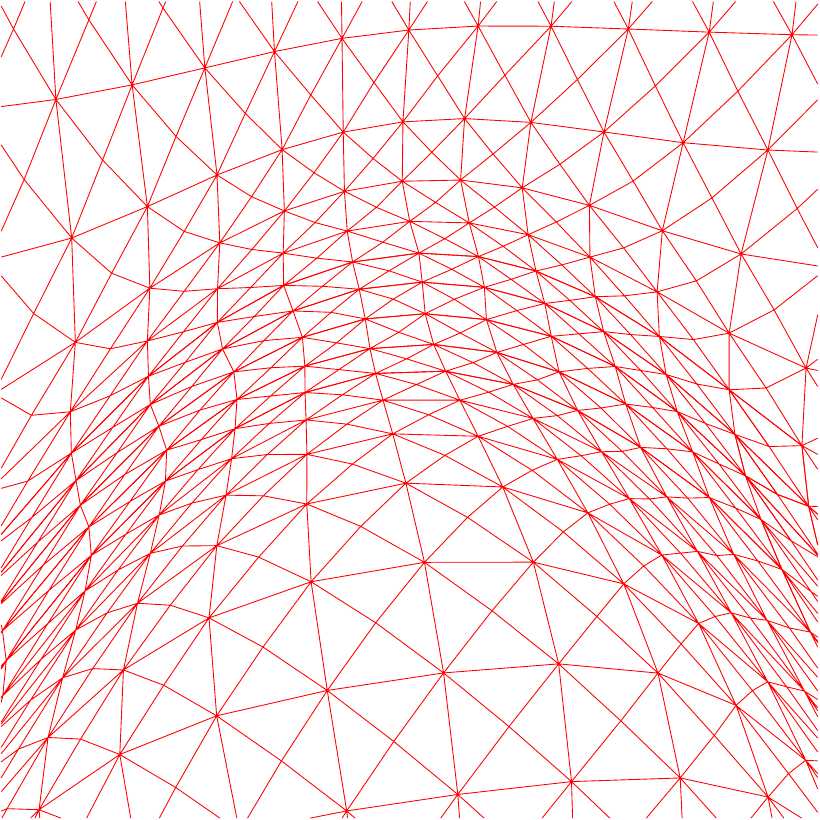}%
      \hfill%
      \includegraphics[clip, width=0.30\linewidth]{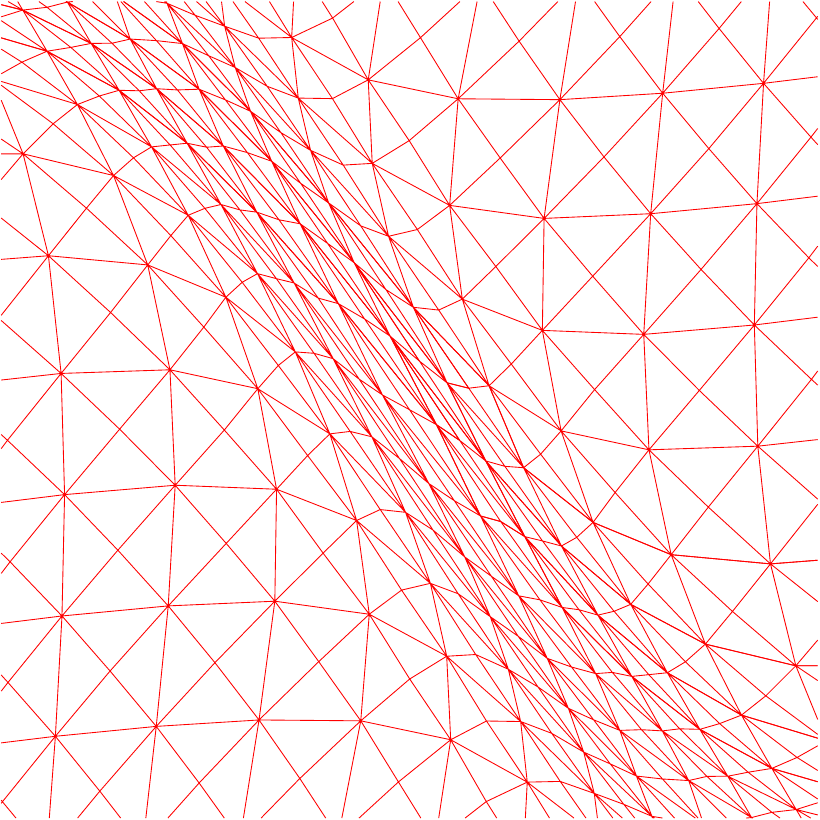}%
      \caption{Huang and Russell's \cref{hr-1}}\label{fig:ex:1:HR}
   \end{subfigure}%
   \\[2.0ex]
   \begin{subfigure}[t]{0.31\linewidth}
      \begin{tikzpicture}
         \begin{semilogxaxis}[
            width=0.75\linewidth,
            height=0.60\linewidth,
            ymin = 0.9, ymax = 2.8, 
         ]
            \addplot[color=\cWinslow, \ltWinslow, mark = +,]
            table [x index=1, y index=4, col sep = space] {results-1-winslow.dat};
            \addlegendentry{Winslow}
            \addplot[color=\cHuang, \ltHuang, mark = +,]
            table [x index=1, y index=4, col sep = space] {results-1-huang.dat};
            \addlegendentry{Huang}
            \addplot[color=\cHR, \ltHR, mark = +,]
            table [x index=1, y index=4, col sep = space] {results-1-hr.dat};
            \addlegendentry{HR}
         \end{semilogxaxis}
      \end{tikzpicture}
      \caption{$Q_{eq}$ vs. $N$}\label{fig:ex:1:qeq}
   \end{subfigure}%
   \hfill%
   \begin{subfigure}[t]{0.31\linewidth}
        \begin{tikzpicture}
         \begin{semilogxaxis}[
            width=0.75\linewidth,
            height=0.60\linewidth,
            ymin = 0.9, ymax = 2.8, 
         ]
            \addplot[color=\cWinslow, \ltWinslow, mark = +,]
            table [x index=1, y index=5, col sep = space] {results-1-winslow.dat};
            \addlegendentry{Winslow}
            \addplot[color=\cHuang, \ltHuang, mark = +,]
            table [x index=1, y index=5, col sep = space] {results-1-huang.dat};
            \addlegendentry{Huang}
            \addplot[color=\cHR, \ltHR, mark = +,]
            table [x index=1, y index=5, col sep = space] {results-1-hr.dat};
            \addlegendentry{HR}
         \end{semilogxaxis}
      \end{tikzpicture}
      \caption{$Q_{ali}$ vs. $N$}\label{fig:ex:1:qali}
   \end{subfigure}%
   \hfill%
   \begin{subfigure}[t]{0.31\linewidth}
      \begin{tikzpicture}
         \begin{loglogaxis}[
            width=0.75\linewidth,
            height=0.60\linewidth,
            legend style={at={(0.98,1.00)}, anchor=north east,align=left},
         ]
            \addplot[color=\cWinslow, \ltWinslow, mark = +,]
            table [x index=1, y index=2, col sep = space] {results-1-winslow.dat};
            \addlegendentry{Winslow}
            \addplot[color=\cHuang, \ltHuang, mark = +,]
            table [x index=1, y index=2, col sep = space] {results-1-huang.dat};
            \addlegendentry{Huang}
            \addplot[color=\cHR, \ltHR, mark = +,]
            table [x index=1, y index=2, col sep = space] {results-1-hr.dat};
            \addlegendentry{HR}
            \addplot[color=\cN, \ltN, mark = none] coordinates
            {(1000, 0.005) (10000, 0.0005) (50000, 0.0001)};
        \end{loglogaxis}
      \end{tikzpicture}
      \caption{$L^2$ interpolation error vs. $N$}\label{fig:ex:1:error}
   \end{subfigure}%
   \caption{\Cref{ex:1}: example meshes (left),
   close-ups near the wave tip (middle) and in the middle (right),
   mesh quality measures, and $L^2$ interpolation error (black line represents $N^{-1}$).}\label{fig:ex:1}%
\end{figure}
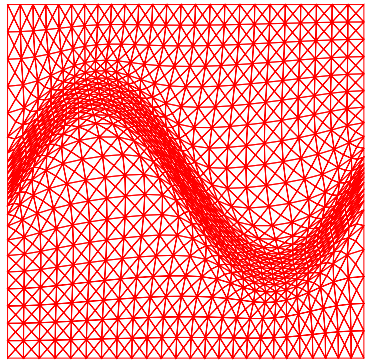
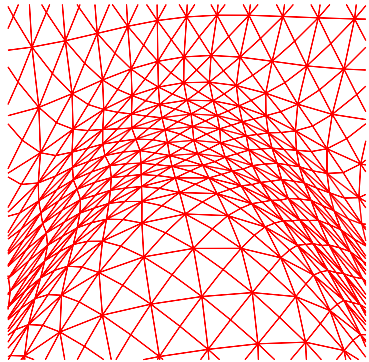
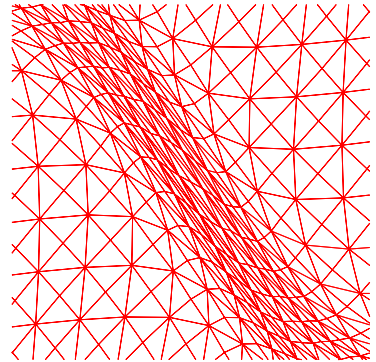
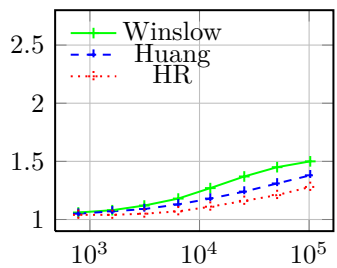
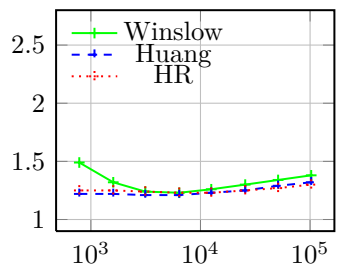
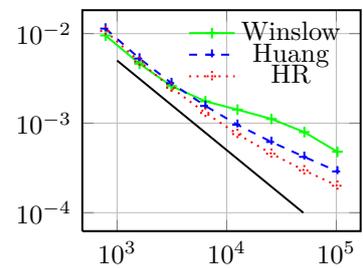

\begin{figure}[p]
   \begin{subfigure}[t]{1.0\linewidth}
      \includegraphics[clip, width=0.30\linewidth]{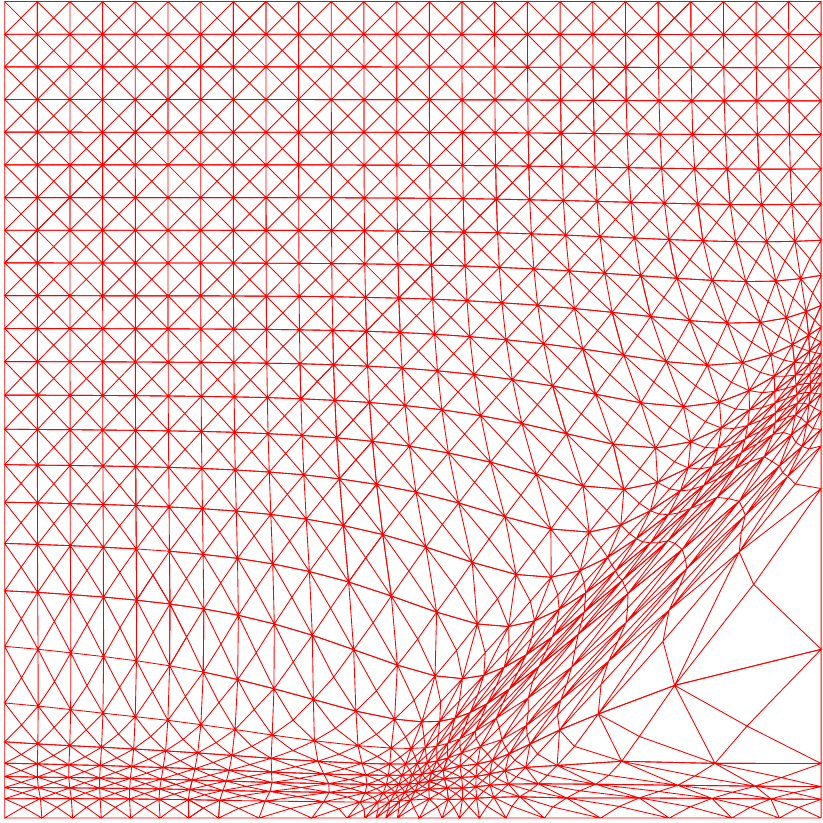}%
      \hfill%
      \includegraphics[clip, width=0.30\linewidth]{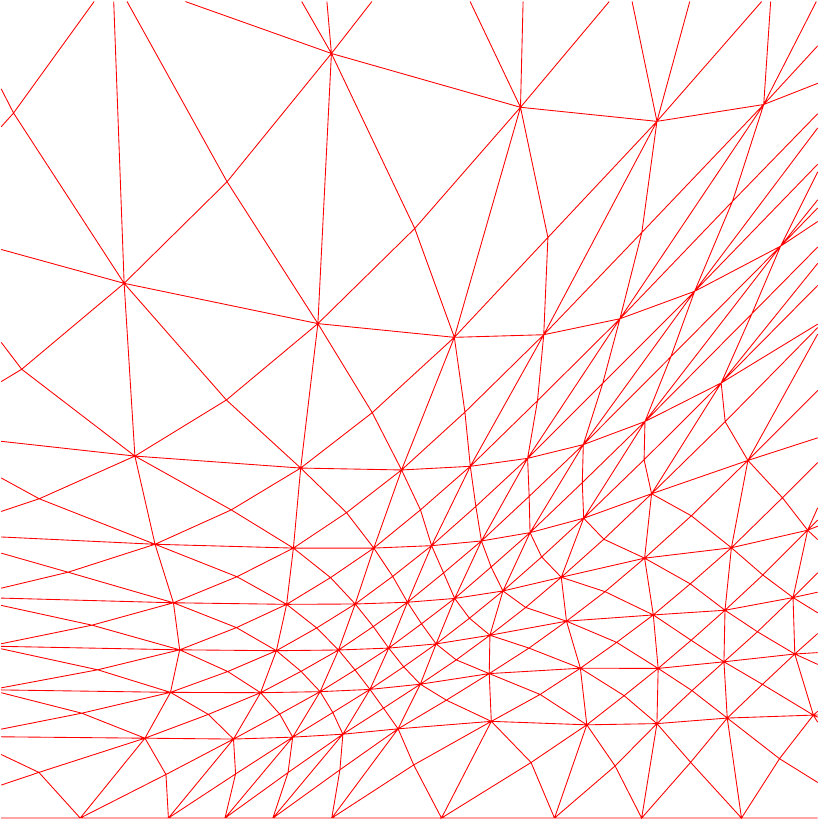}%
      \hfill%
      \includegraphics[clip, width=0.30\linewidth]{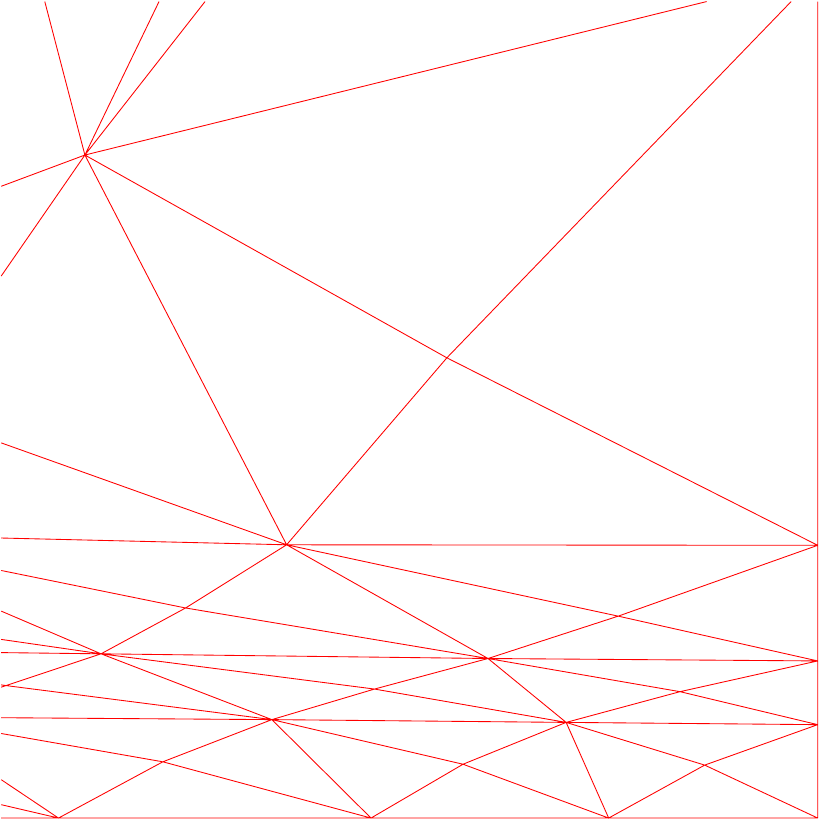}%
      \caption{Winslow's \cref{win-2}}\label{fig:ex:2:winslow}
   \end{subfigure}%
   \\[2.0ex]
   \begin{subfigure}[t]{1.0\linewidth}
      \includegraphics[clip, width=0.30\linewidth]{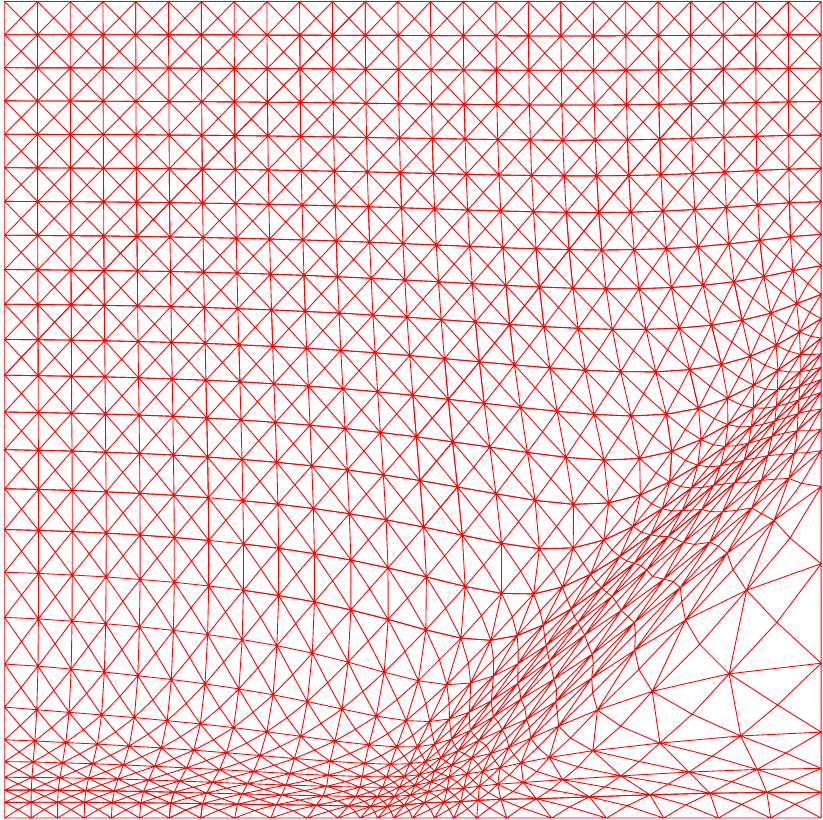}%
      \hfill%
      \includegraphics[clip, width=0.30\linewidth]{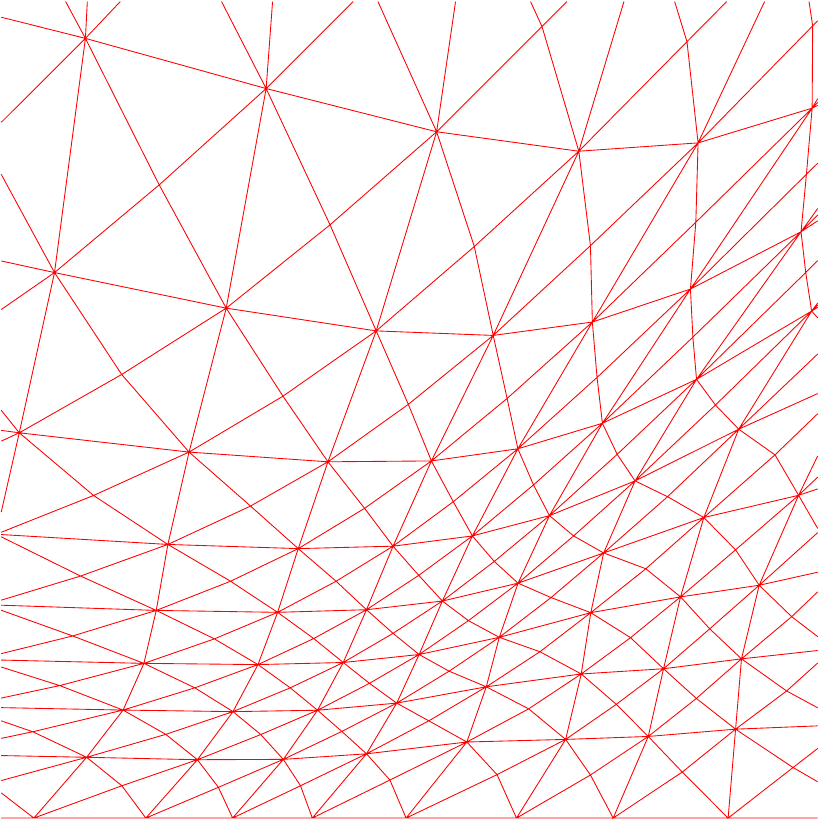}%
      \hfill%
      \includegraphics[clip, width=0.30\linewidth]{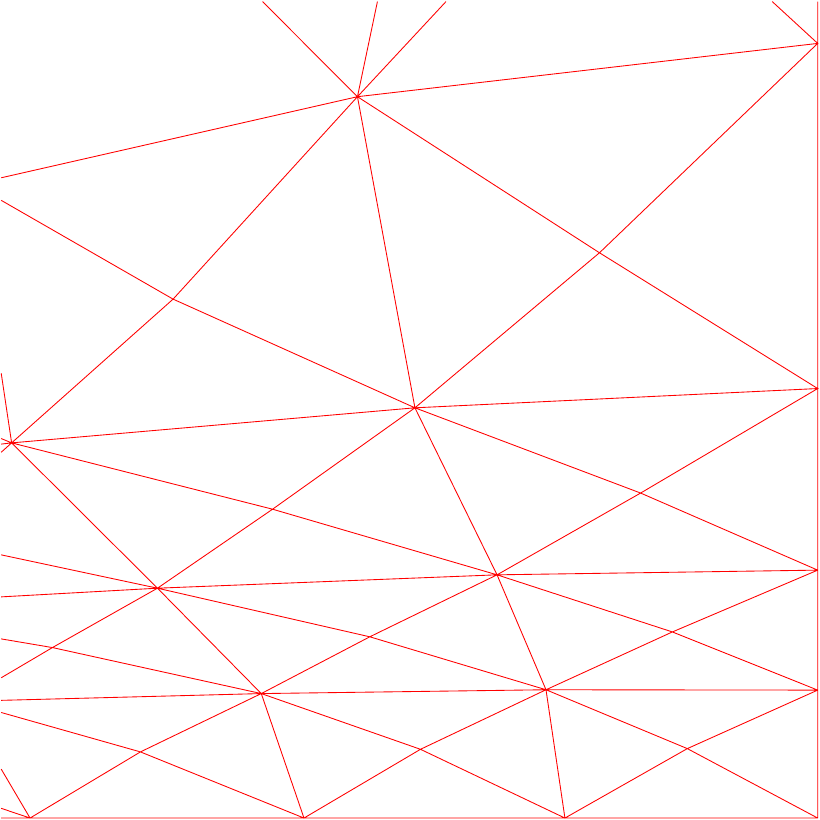}%
      \caption{using Huang's \cref{huang-1}}\label{fig:ex:2:huang}
   \end{subfigure}%
   \\[2.0ex]
   \begin{subfigure}[t]{1.0\linewidth}
      \includegraphics[clip, width=0.30\linewidth]{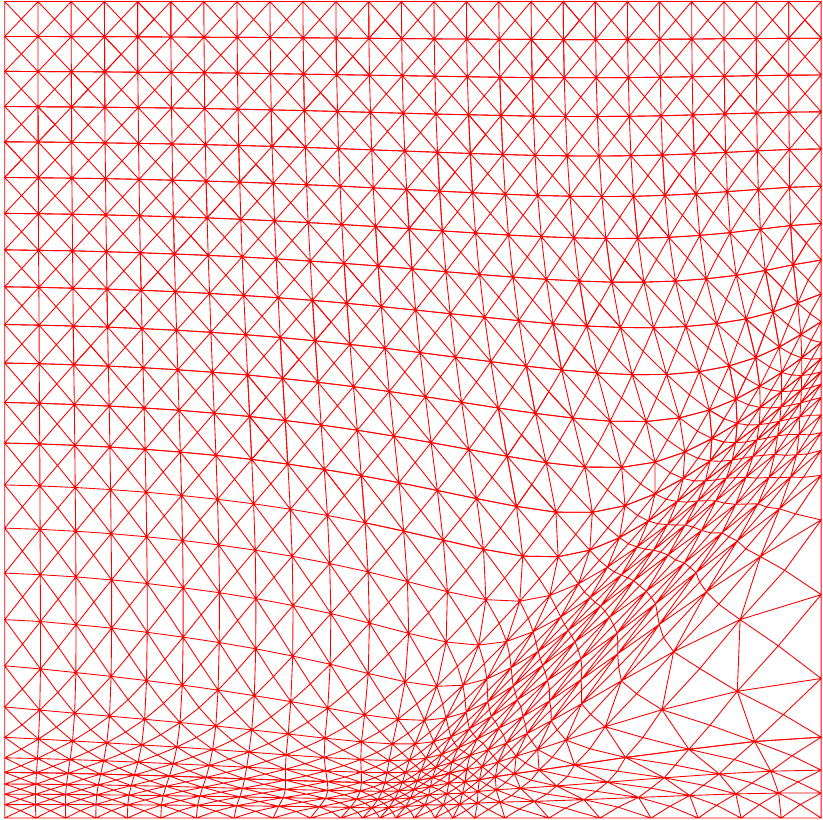}%
      \hfill%
      \includegraphics[clip, width=0.30\linewidth]{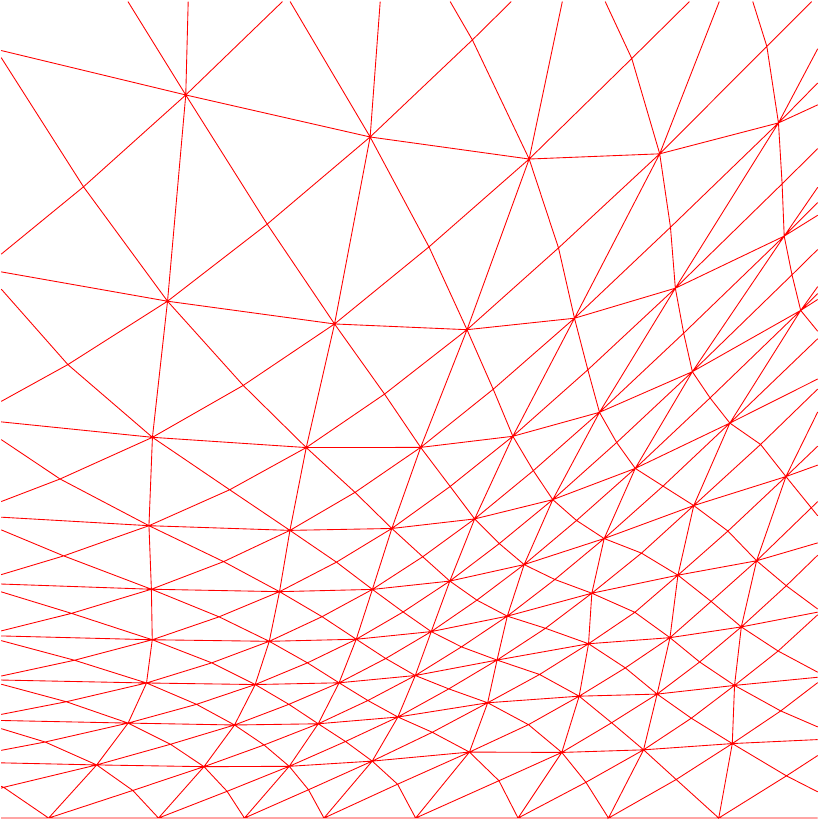}%
      \hfill%
      \includegraphics[clip, width=0.30\linewidth]{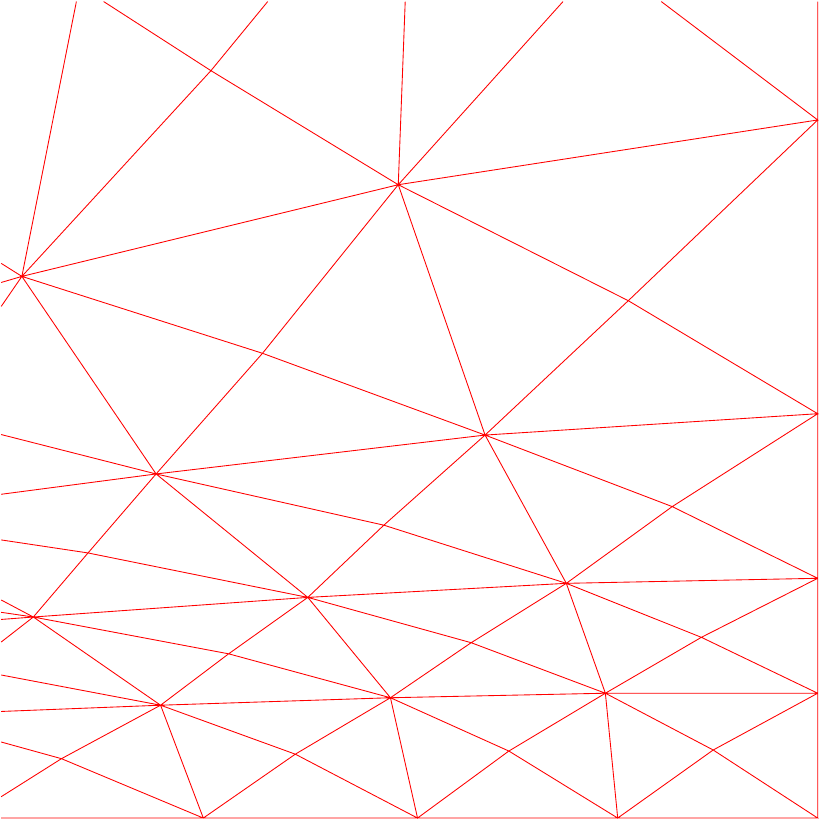}%
      \caption{Huang and Russell's \cref{hr-1}}\label{fig:ex:2:HR}
   \end{subfigure}%
   \\[2.0ex]
   \begin{subfigure}[t]{0.31\linewidth}
      \begin{tikzpicture}
         \begin{semilogxaxis}[
            width=0.75\linewidth,
            height=0.60\linewidth,
            ymin = 0.9, ymax = 2.8, 
         ]
            \addplot[color=\cWinslow, \ltWinslow, mark = +,]
            table [x index=1, y index=4, col sep = space] {results-2-winslow.dat};
            \addlegendentry{Winslow}
            \addplot[color=\cHuang, \ltHuang, mark = +,]
            table [x index=1, y index=4, col sep = space] {results-2-huang.dat};
            \addlegendentry{Huang}
            \addplot[color=\cHR, \ltHR, mark = +,]
            table [x index=1, y index=4, col sep = space] {results-2-hr.dat};
            \addlegendentry{HR}
         \end{semilogxaxis}
      \end{tikzpicture}
      \caption{$Q_{eq}$ vs. $N$}\label{fig:ex:2:qeq}
   \end{subfigure}%
   \hfill%
   \begin{subfigure}[t]{0.31\linewidth}
        \begin{tikzpicture}
         \begin{semilogxaxis}[
            width=0.75\linewidth,
            height=0.60\linewidth,
            ymin = 0.9, ymax = 2.8, 
         ]
            \addplot[color=\cWinslow, \ltWinslow, mark = +,]
            table [x index=1, y index=5, col sep = space] {results-2-winslow.dat};
            \addlegendentry{Winslow}
            \addplot[color=\cHuang, \ltHuang, mark = +,]
            table [x index=1, y index=5, col sep = space] {results-2-huang.dat};
            \addlegendentry{Huang}
            \addplot[color=\cHR, \ltHR, mark = +,]
            table [x index=1, y index=5, col sep = space] {results-2-hr.dat};
            \addlegendentry{HR}
         \end{semilogxaxis}
      \end{tikzpicture}
      \caption{$Q_{ali}$ vs. $N$}\label{fig:ex:2:qali}
   \end{subfigure}%
   \hfill%
   \begin{subfigure}[t]{0.31\linewidth}
      \begin{tikzpicture}
         \begin{loglogaxis}[
            width=0.75\linewidth,
            height=0.60\linewidth,
            legend style={at={(0.98,1.00)}, anchor=north east,align=left},
         ]
            \addplot[color=\cWinslow, \ltWinslow, mark = +]
            table [x index=1, y index=2, col sep = space] {results-2-winslow.dat};
            \addlegendentry{Winslow}
            \addplot[color=\cHuang, \ltHuang, mark = +,]
            table [x index=1, y index=2, col sep = space] {results-2-huang.dat};
            \addlegendentry{Huang}
            \addplot[color=\cHR, \ltHR, mark = +,]
            table [x index=1, y index=2, col sep = space] {results-2-hr.dat};
            \addlegendentry{HR}
            \addplot[color=\cN, \ltN, mark = none] coordinates
            {(1000, 0.004) (10000, 0.0004) (40000, 0.0001)};
        \end{loglogaxis}
      \end{tikzpicture}
      \caption{$L^2$ error vs. $N$}\label{fig:ex:2:error}
   \end{subfigure}%
   \caption{\Cref{ex:2}: example meshes (left),
   close-ups near the wave meeting the boundary layer (middle)
   and in the right bottom corner (right),
   mesh quality measures, and $L^2$ interpolation error  (black line represents $N^{-1}$).}\label{fig:ex:2}%
\end{figure}
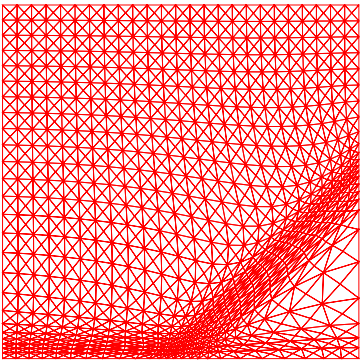
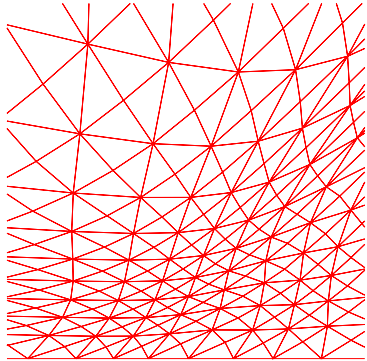
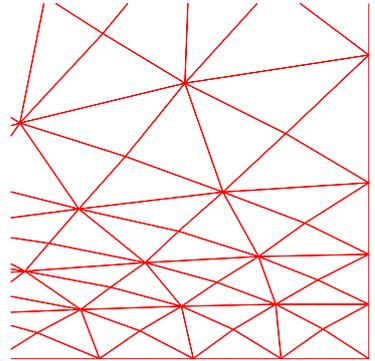
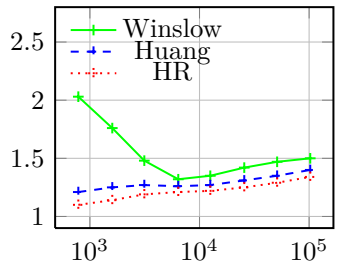
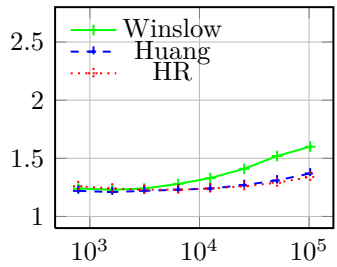
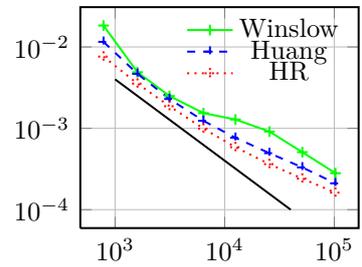

First, we consider two dimensional meshes constructed for the unit square $\Omega = (0,1) \times (0,1)$ and the test functions
\begin{example}
\label{ex:1}
\[
   u = \tanh
      \left(-30 \left[ y - 0.5 - 0.25 \sin (2 \pi x) \right] \right)
   .
\]
\end{example}
\begin{example}
\label{ex:2}
\[
   u = \tanh(25 y) - \tanh(25 ( x - y - 0.5))
   .
\]
\end{example}
Example meshes, close-ups, as well as the mesh quality measures and the $L^2$ interpolation error are given in \cref{fig:ex:1,fig:ex:2}.

For these examples, all three functionals provide good size and shape adaptation.
A closer look at the mesh quality measures shows that, although all three functionals provide comparable meshes which are reasonably close to the prescribed metric tensor, meshes constructed using H and HR functionals have better correspondence to the prescribed metric tensor.
In both two-dimensional examples, H and HR functionals provide very similar grids which are closer to the prescribed size and shape (that is, smaller values of $Q_{eq}$ and $Q_{ali}$).
This is also reflected in the error of the linear interpolation: HR functional \cref{hr-1} provides the smallest error, followed by H functional \cref{huang-1} and then W functional \cref{win-2}.
It seems that W functional is a bit too aggressive in moving nodes towards the neighborhood of areas of interest, providing a higher density of the nodes along the anisotropic features of the given function while coarsening out the mesh nearby, leading to a steeper element size gradation.
Interestingly, for both examples the convergence of the linear interpolation error for W functional (\cref{fig:ex:1:error,fig:ex:2:error}) slows down near $N = 10^4$ and returns to the order $\mathcal{O}(N^{-1}) $ as the mesh is refined ($N$ is the number of mesh elements).
It is unclear to us what causes this for W functional.

\begin{figure}[p]
   \begin{subfigure}[t]{1.0\linewidth}
      \includegraphics[clip, width=0.30\linewidth]{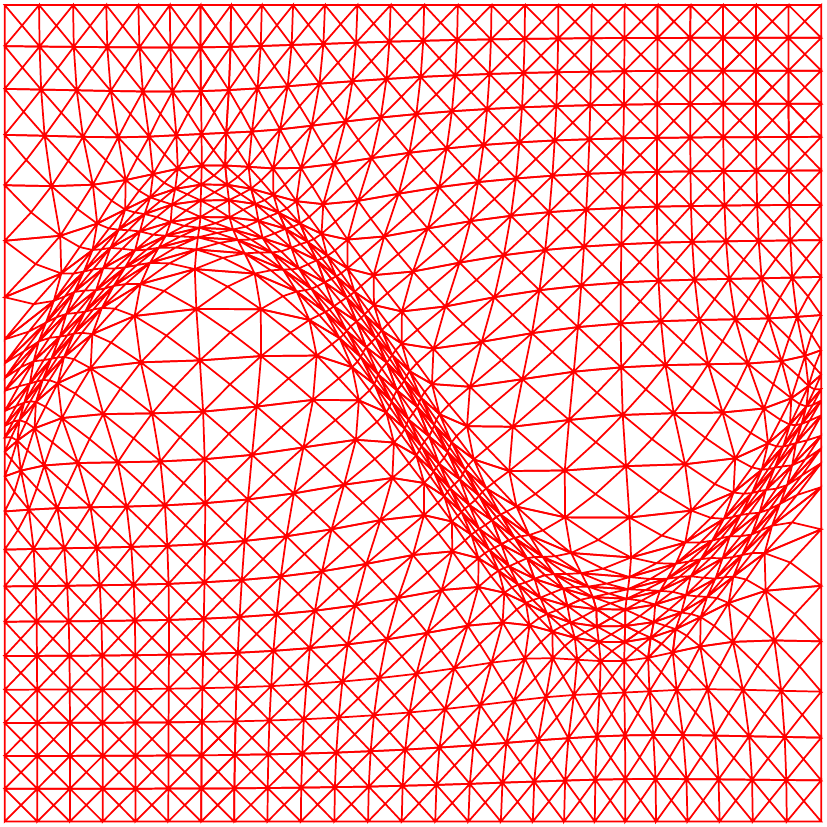}%
      \hfill%
      \includegraphics[clip, width=0.30\linewidth]{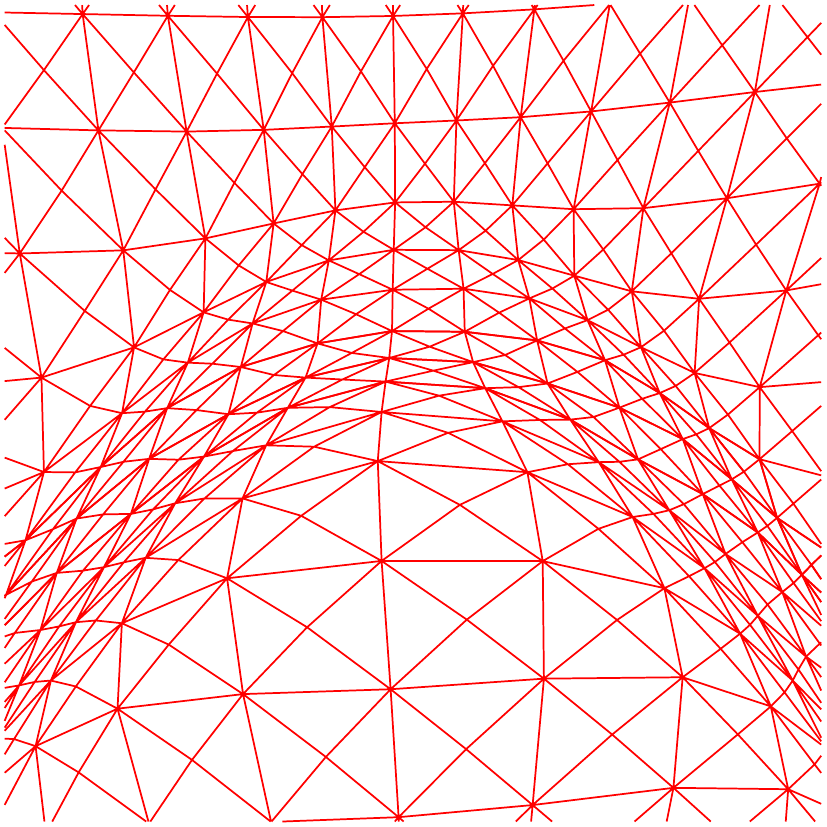}%
      \hfill%
      \includegraphics[clip, width=0.30\linewidth]{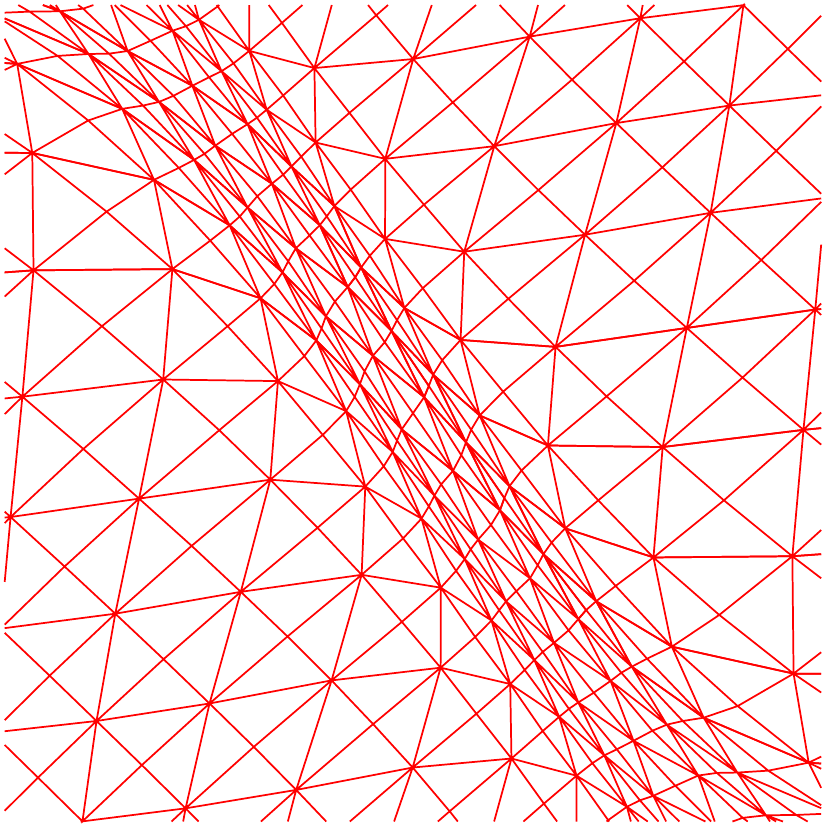}%
      \caption{Winslow's \cref{win-2}}\label{fig:ex:1:h1:winslow}
   \end{subfigure}%
   \\[2.0ex]
   \begin{subfigure}[t]{1.0\linewidth}
      \includegraphics[clip, width=0.30\linewidth]{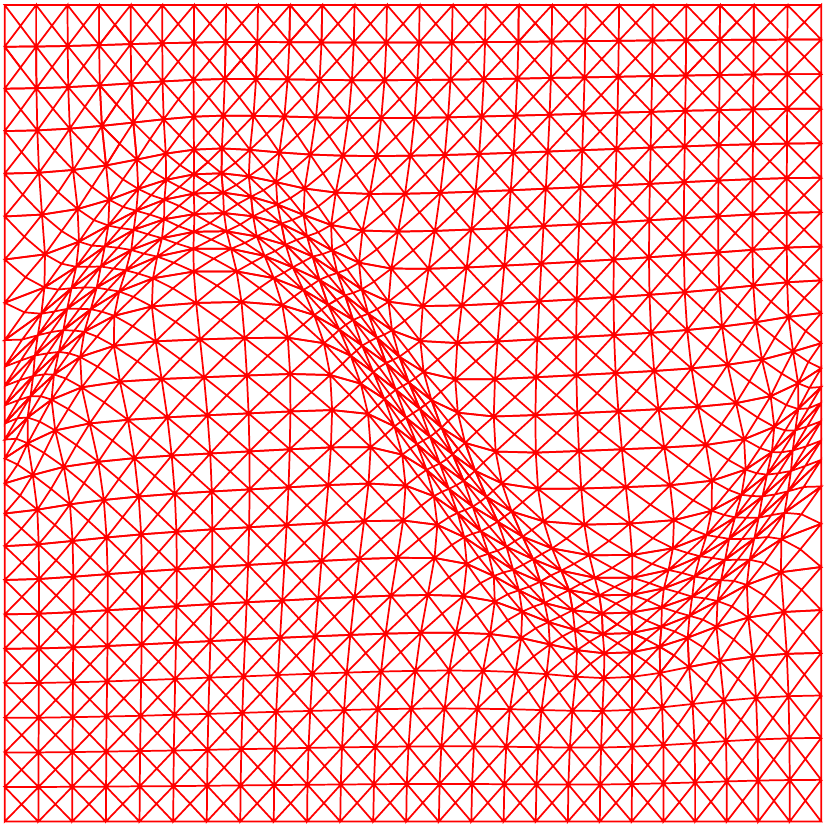}%
      \hfill%
      \includegraphics[clip, width=0.30\linewidth]{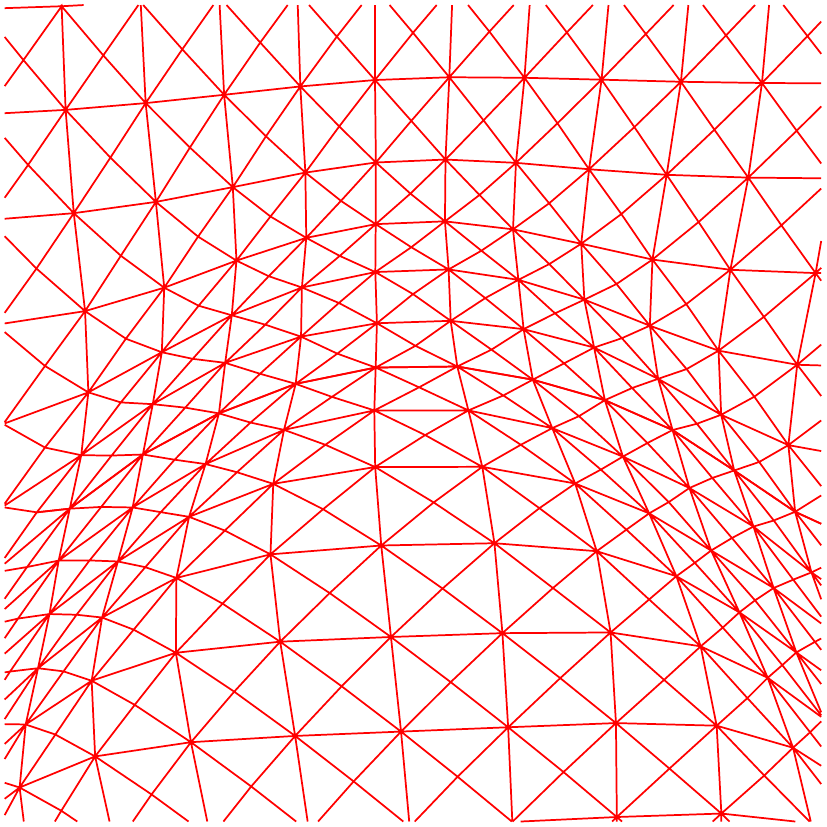}%
      \hfill%
      \includegraphics[clip, width=0.30\linewidth]{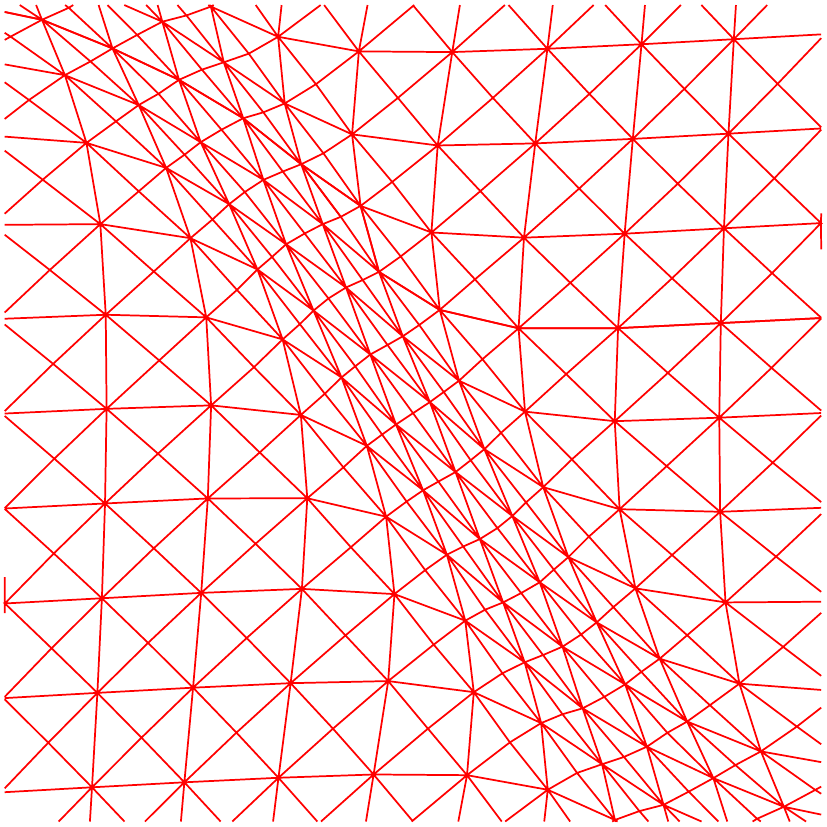}%
      \caption{Huang's \cref{huang-1}}\label{fig:ex:1:h1:huang}
   \end{subfigure}%
   \\[2.0ex]
   \begin{subfigure}[t]{1.0\linewidth}
      \includegraphics[clip, width=0.30\linewidth]{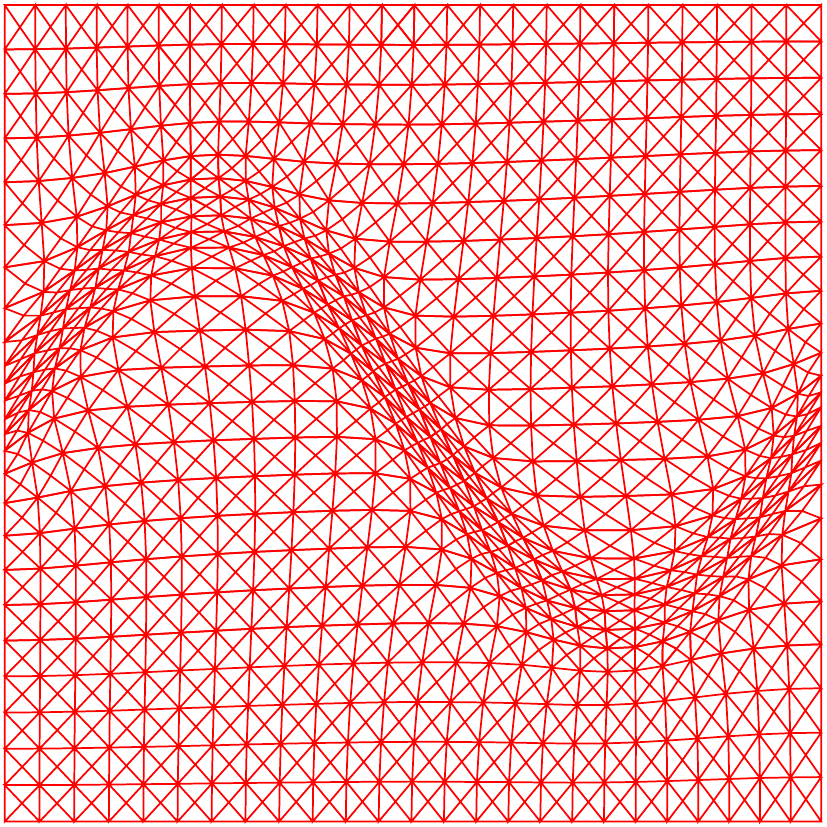}%
      \hfill%
      \includegraphics[clip, width=0.30\linewidth]{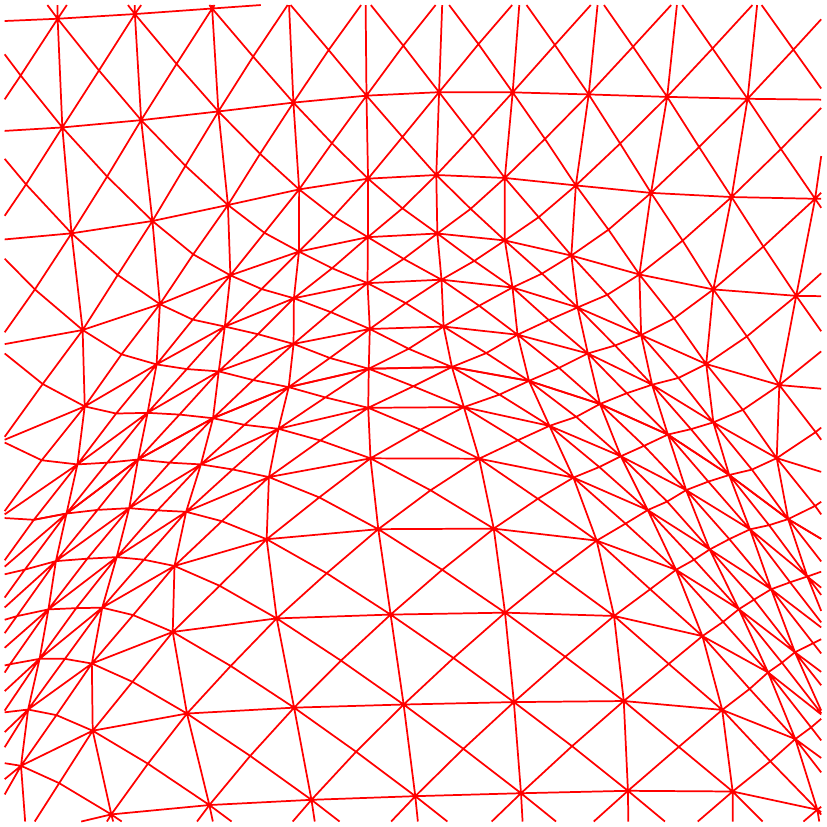}%
      \hfill%
      \includegraphics[clip, width=0.30\linewidth]{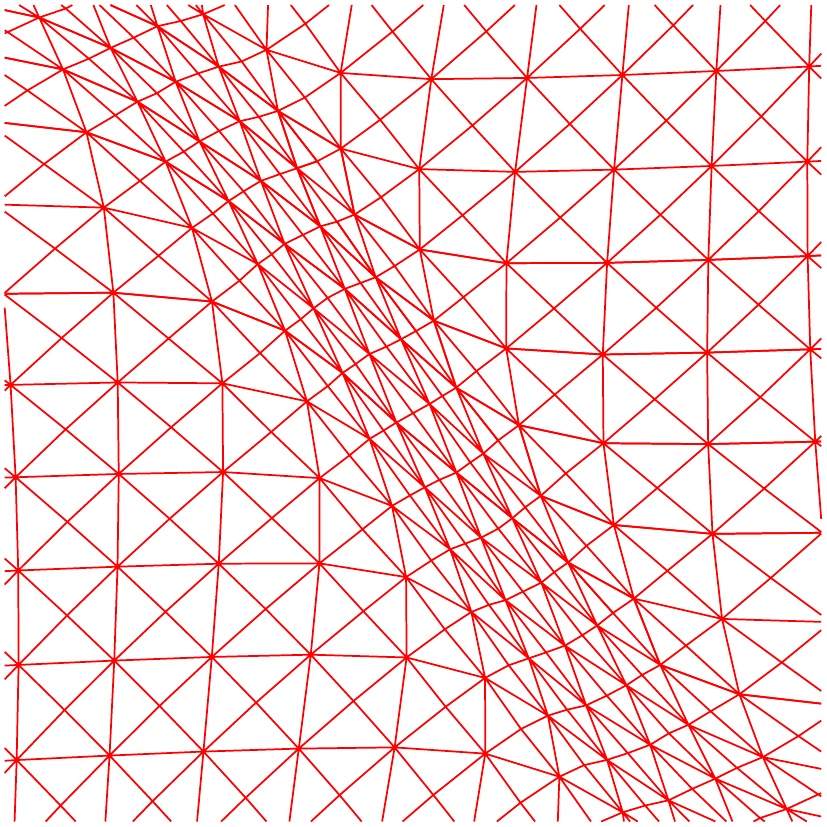}%
      \caption{Huang and Russell's \cref{hr-1}}\label{fig:ex:1:h1:HR}
   \end{subfigure}%
   \\[2.0ex]
   \begin{subfigure}[t]{0.31\linewidth}
      \begin{tikzpicture}
         \begin{semilogxaxis}[
            width=0.75\linewidth,
            height=0.60\linewidth,
            ymin = 0.9, ymax = 2.8, 
         ]
            \addplot[color=\cWinslow, \ltWinslow, mark = +,]
            table [x index=1, y index=4, col sep = space] {results-1-winslow-h1.dat};
            \addlegendentry{Winslow}
            \addplot[color=\cHuang, \ltHuang, mark = +,]
            table [x index=1, y index=4, col sep = space] {results-1-huang-h1.dat};
            \addlegendentry{Huang}
            \addplot[color=\cHR, \ltHR, mark = +,]
            table [x index=1, y index=4, col sep = space] {results-1-hr-h1.dat};
            \addlegendentry{HR}
         \end{semilogxaxis}
      \end{tikzpicture}
      \caption{$Q_{eq}$ vs. $N$}\label{fig:ex:1:h1:qeq}
   \end{subfigure}%
   \hfill%
   \begin{subfigure}[t]{0.31\linewidth}
        \begin{tikzpicture}
         \begin{semilogxaxis}[
            width=0.75\linewidth,
            height=0.60\linewidth,
            ymin = 0.9, ymax = 2.8, 
         ]
            \addplot[color=\cWinslow, \ltWinslow, mark = +,]
            table [x index=1, y index=5, col sep = space] {results-1-winslow-h1.dat};
            \addlegendentry{Winslow}
            \addplot[color=\cHuang, \ltHuang, mark = +,]
            table [x index=1, y index=5, col sep = space] {results-1-huang-h1.dat};
            \addlegendentry{Huang}
            \addplot[color=\cHR, \ltHR, mark = +,]
            table [x index=1, y index=5, col sep = space] {results-1-hr-h1.dat};
            \addlegendentry{HR}
         \end{semilogxaxis}
      \end{tikzpicture}
      \caption{$Q_{ali}$ vs. $N$}\label{fig:ex:1:h1:qali}
   \end{subfigure}%
   \hfill%
   \begin{subfigure}[t]{0.31\linewidth}
      \begin{tikzpicture}
         \begin{loglogaxis}[
            width=0.75\linewidth,
            height=0.60\linewidth,
            legend style={at={(0.98,1.00)}, anchor=north east,align=left},
         ]
            \addplot[color=\cWinslow, \ltWinslow, mark = +,]
            table [x index=1, y index=2, col sep = space] {results-1-winslow-h1.dat};
            \addlegendentry{Winslow}
            \addplot[color=\cHuang, \ltHuang, mark = +,]
            table [x index=1, y index=2, col sep = space] {results-1-huang-h1.dat};
            \addlegendentry{Huang}
            \addplot[color=\cHR, \ltHR, mark = +,]
            table [x index=1, y index=2, col sep = space] {results-1-hr-h1.dat};
            \addlegendentry{HR}
            \addplot[color=\cN, \ltN, mark = none] coordinates
            {(1000, 0.005) (10000, 0.0005) (50000, 0.0001)};
        \end{loglogaxis}
      \end{tikzpicture}
      \caption{$L^2$ interpolation error vs. $N$}\label{fig:ex:1:h1:error}
   \end{subfigure}%
   \caption{\Cref{ex:1} using a metric tensor for the $H^1$ semi-norm:
   example meshes (left),
   close-ups near the wave tip (middle) and in the middle (right),
   mesh quality measures, and $L^2$ interpolation error
   (black line represents $N^{-1}$).}\label{fig:ex:1:h1}%
\end{figure}

For \cref{ex:1} we also compute adaptive meshes for the metric tensor which is optimal for the $H^1$ semi norm error (see, e.g.,\ \cite{Hua06} for details on the metric tensor).
The results shown in \cref{fig:ex:1:h1} are essentially the same as in \cref{fig:ex:1}: HR functional \cref{hr-1} provides the smallest error, followed by H functional \cref{huang-1} and then W functional \cref{win-2}.
For this metric tensor, the $L^2$ error for the H and HR functionals (\cref{fig:ex:1:h1:error}) is slightly larger than in \cref{fig:ex:1:error}, which is not surprising, since the metric tensor chosen for \cref{fig:ex:1} is optimal for the $L^2$ error.
Thus, adding the equidistribution property to the meshing functional seems to help to obtain meshes which are closer to fulfilling the prescribed properties.
Interestingly, the $L^2$ error for the W functional seems to be a bit smaller if the $H^1$ metric tensor is used.
This can be explained by the fact that the W functional does not incorporate the equidistribution property and, thus, doesn't exactly generate a mesh which is uniform with respect to the provided metric tensor.
Thus, it is not quite clear if one is able to generate the optimal mesh when using the W functional.

\subsection{Three dimensions}
In three dimensions, we consider the unit cube $\Omega = (0,1)\times (0,1)\times (0,1)$ and the following test functions.

\begin{example}
\label{ex:3}
\begin{align*}
   u= &\tanh \left(30\left[ (4x-2.0)^2 +(4y-2.0)^2 +(4z-2.0)^2 -0.1875 \right] \right)\\
    + &\tanh \left(30\left[ (4x-2.5)^2 +(4y-2.5)^2 +(4z-2.5)^2 -0.1875 \right] \right)\\
    + &\tanh \left(30\left[ (4x-2.5)^2 +(4y-1.5)^2 +(4z-2.5)^2 -0.1875 \right] \right)\\
    + &\tanh \left(30\left[ (4x-1.5)^2 +(4y-2.5)^2 +(4z-2.5)^2 -0.1875 \right] \right)\\
    + &\tanh \left(30\left[ (4x-1.5)^2 +(4y-1.5)^2 +(4z-2.5)^2 -0.1875 \right] \right)\\
    + &\tanh \left(30\left[ (4x-2.5)^2 +(4y-2.5)^2 +(4z-1.5)^2 -0.1875 \right] \right)\\
    + &\tanh \left(30\left[ (4x-2.5)^2 +(4y-1.5)^2 +(4z-1.5)^2 -0.1875 \right] \right)\\
    + &\tanh \left(30\left[ (4x-1.5)^2 +(4y-2.5)^2 +(4z-1.5)^2 -0.1875 \right] \right)\\
    + &\tanh \left(30\left[ (4x-1.5)^2 +(4y-1.5)^2 +(4z-1.5)^2 -0.1875 \right] \right)
   .
\end{align*}
\end{example}

\begin{example}
\label{ex:4}
\[
   u = \tanh \Big( - 30 
         [ z - 0.5 - 0.25 \sin (2 \pi x) \sin(\pi y) ]
      \Big)
   .
\]
\end{example}

\begin{example}
\label{ex:5}
\[
   u = \tanh \bigg( -30 
      \Big\{ z- \tanh \big(-30 
                  \left [y - 0.5 - 0.25 \sin (2 \pi x)\right ] \big) \Big\} \bigg)
   .
\]
\end{example}

\begin{figure}[p]
   \begin{subfigure}[t]{0.31\linewidth}
      \includegraphics[clip, width=1.0\linewidth]{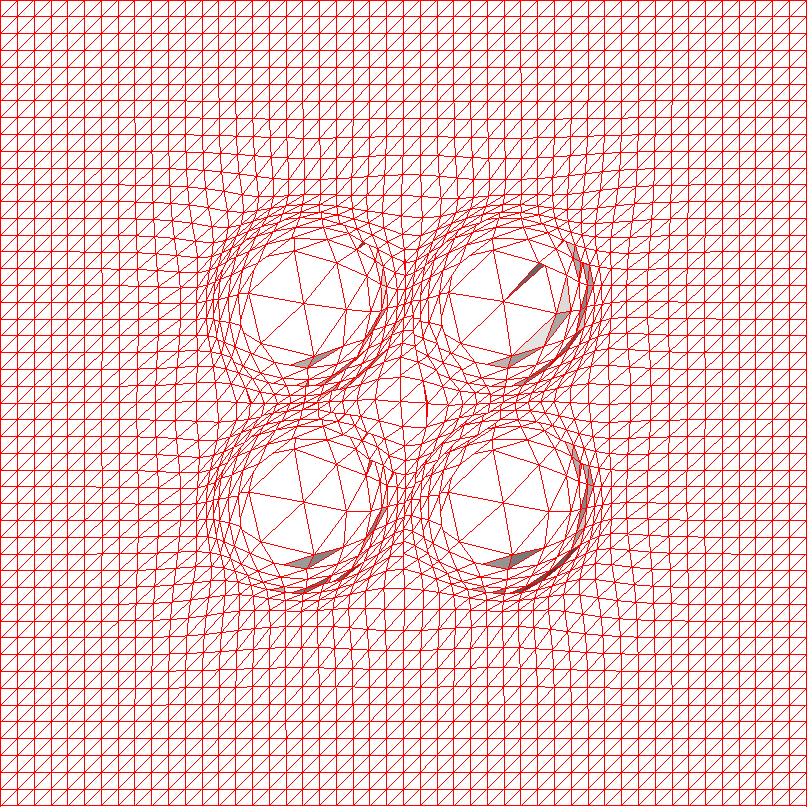}%
      \caption{Winslow's \cref{win-2}}\label{fig:ex:3:winslow}
   \end{subfigure}%
   \hfill%
   \begin{subfigure}[t]{0.31\linewidth}
      \includegraphics[clip, width=1.0\linewidth]{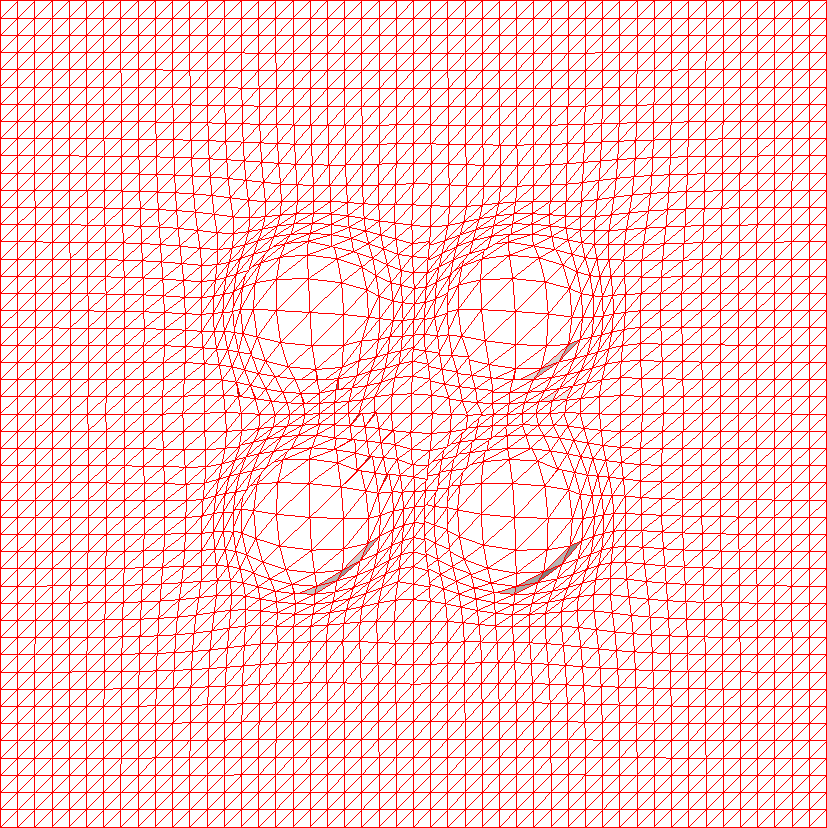}%
      \caption{Huang's \cref{huang-1}}\label{fig:ex:3:huang}
   \end{subfigure}%
   \hfill%
   \begin{subfigure}[t]{0.31\linewidth}
      \includegraphics[clip, width=1.0\linewidth]{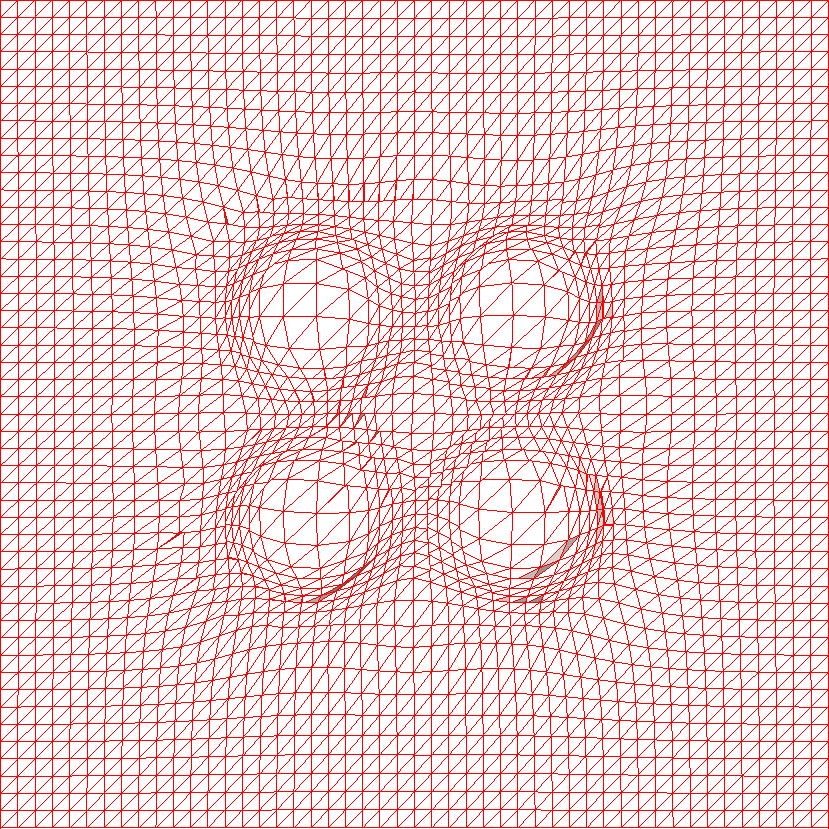}%
      \caption{Huang and Russell's \cref{hr-1}}\label{fig:ex:3:HR}
   \end{subfigure}%
   \\[2.0ex]
   \begin{subfigure}[t]{0.31\linewidth}
      \includegraphics[clip, width=1.0\linewidth]{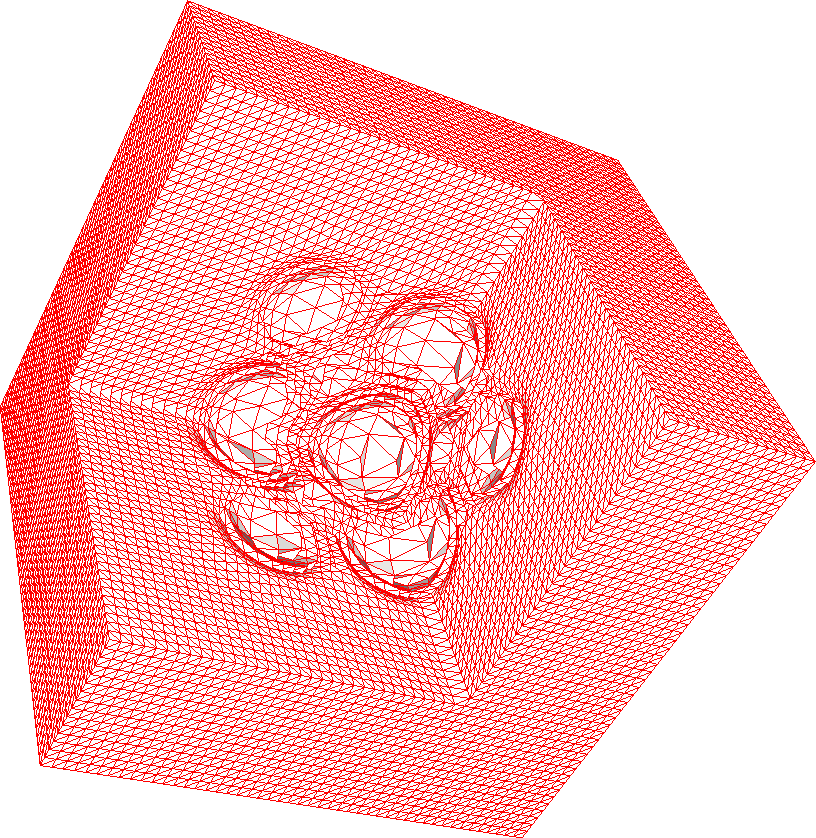}%
      \caption{Winslow's \cref{win-2}}\label{fig:ex:3:winslow+3d}
   \end{subfigure}%
   \hfill%
   \begin{subfigure}[t]{0.31\linewidth}
      \includegraphics[clip, width=1.0\linewidth]{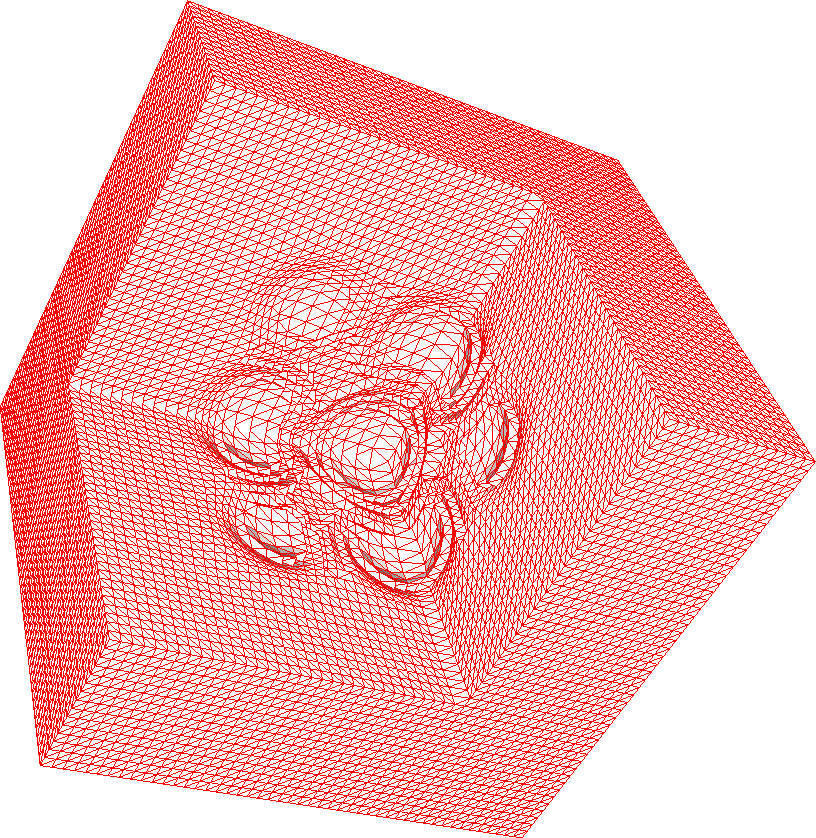}%
      \caption{Huang's \cref{huang-1}}\label{fig:ex:3:huang+3d}
   \end{subfigure}%
   \hfill%
   \begin{subfigure}[t]{0.31\linewidth}
      \includegraphics[clip, width=1.0\linewidth]{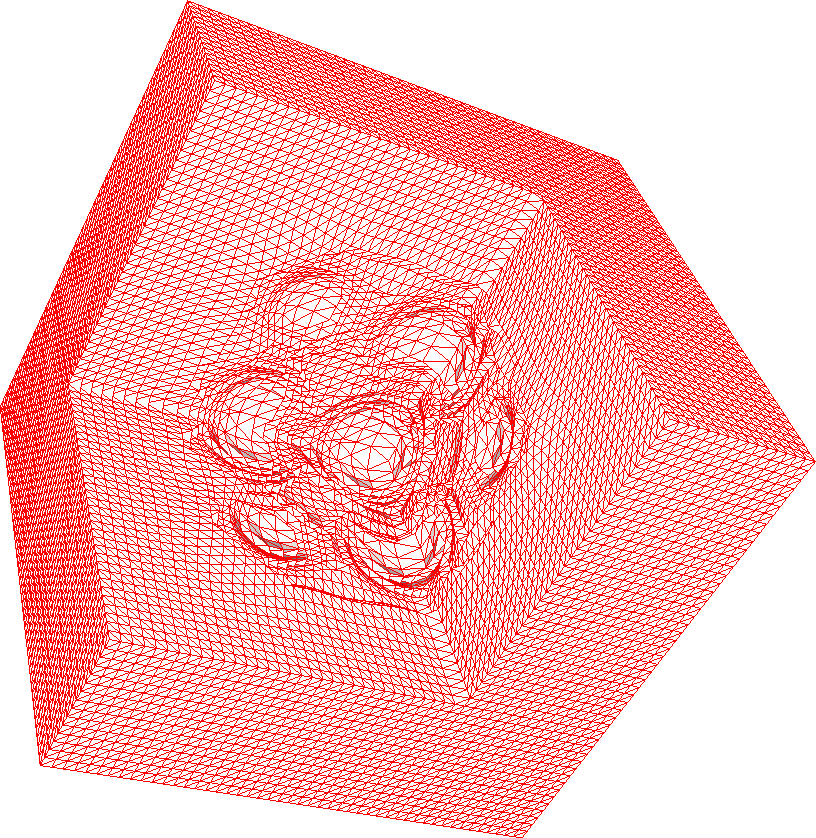}%
      \caption{Huang and Russell's \cref{hr-1}}\label{fig:ex:3:HR+3d}
   \end{subfigure}%
   \\[2.0ex]
   \begin{subfigure}[t]{0.31\linewidth}
      \begin{tikzpicture}
         \begin{semilogxaxis}[
            width=0.75\linewidth,
            height=0.60\linewidth,
            ymin = 0.9, ymax = 2.8, 
         ]
            \addplot[color=\cWinslow, \ltWinslow, mark = +,]
            table [x index=1, y index=4, col sep = space] {results-3-winslow.dat};
            \addlegendentry{Winslow}
            \addplot[color=\cHuang, \ltHuang, mark = +,]
            table [x index=1, y index=4, col sep = space] {results-3-huang.dat};
            \addlegendentry{Huang}
            \addplot[color=\cHR, \ltHR, mark = +,]
            table [x index=1, y index=4, col sep = space] {results-3-hr.dat};
            \addlegendentry{HR}
         \end{semilogxaxis}
      \end{tikzpicture}
      \caption{$Q_{eq}$ vs. $N$}\label{fig:ex:3:qeq}
   \end{subfigure}%
   \hfill%
   \begin{subfigure}[t]{0.31\linewidth}
        \begin{tikzpicture}
         \begin{semilogxaxis}[
            width=0.75\linewidth,
            height=0.60\linewidth,
            ymin = 0.9, ymax = 2.8, 
         ]
            \addplot[color=\cWinslow, \ltWinslow, mark = +,]
            table [x index=1, y index=5, col sep = space] {results-3-winslow.dat};
            \addlegendentry{Winslow}
            \addplot[color=\cHuang, \ltHuang, mark = +,]
            table [x index=1, y index=5, col sep = space] {results-3-huang.dat};
            \addlegendentry{Huang}
            \addplot[color=\cHR, \ltHR, mark = +,]
            table [x index=1, y index=5, col sep = space] {results-3-hr.dat};
            \addlegendentry{HR}
         \end{semilogxaxis}
      \end{tikzpicture}
      \caption{$Q_{ali}$ vs. $N$}\label{fig:ex:3:qali}
   \end{subfigure}%
   \hfill%
   \begin{subfigure}[t]{0.31\linewidth}
      \begin{tikzpicture}
         \begin{loglogaxis}[
            width=0.75\linewidth,
            height=0.60\linewidth,
            legend style={at={(0.98,1.00)}, anchor=north east,align=left},
         ]
            \addplot[color=\cWinslow, \ltWinslow, mark = +,]
            table [x index=1, y index=2, col sep = space] {results-3-winslow.dat};
            \addlegendentry{Winslow}
            \addplot[color=\cHuang, \ltHuang, mark = +,]
            table [x index=1, y index=2, col sep = space] {results-3-huang.dat};
            \addlegendentry{Huang}
            \addplot[color=\cHR, \ltHR, mark = +,]
            table [x index=1, y index=2, col sep = space] {results-3-hr.dat};
            \addlegendentry{HR}
            \addplot[color=\cN, \ltN, mark = none] coordinates
            {(10000, 0.109) (40000, 0.043) (160000, 0.017)};
        \end{loglogaxis}
      \end{tikzpicture}
      \caption{$L^2$ error vs. $N$}\label{fig:ex:3:error}
   \end{subfigure}%
   \caption{\Cref{ex:3}. The top row: slice cuts of the meshes.
   The middle row: clip cuts of the meshes.
   The bottom row: $Q_{eq}$, $Q_{ali}$,
   and the $L^2$ norm of the linear interpolation error
   (black line represents $N^{-2/3}$).}\label{fig:ex:3}%
\end{figure}
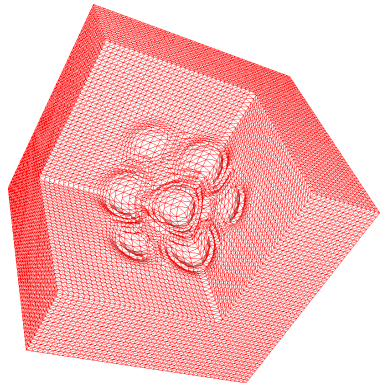
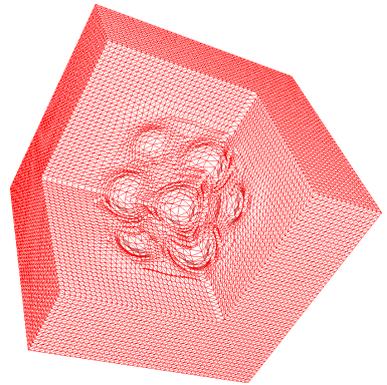
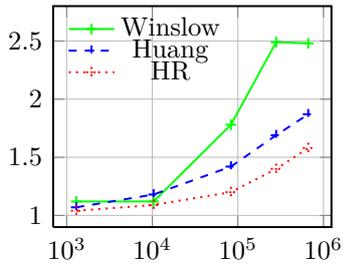
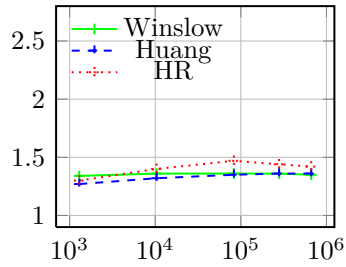
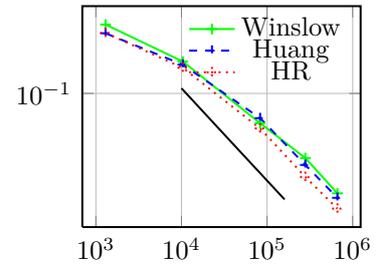

\begin{figure}[p]
   \begin{subfigure}[t]{0.31\linewidth}
      \includegraphics[clip, width=1.0\linewidth]{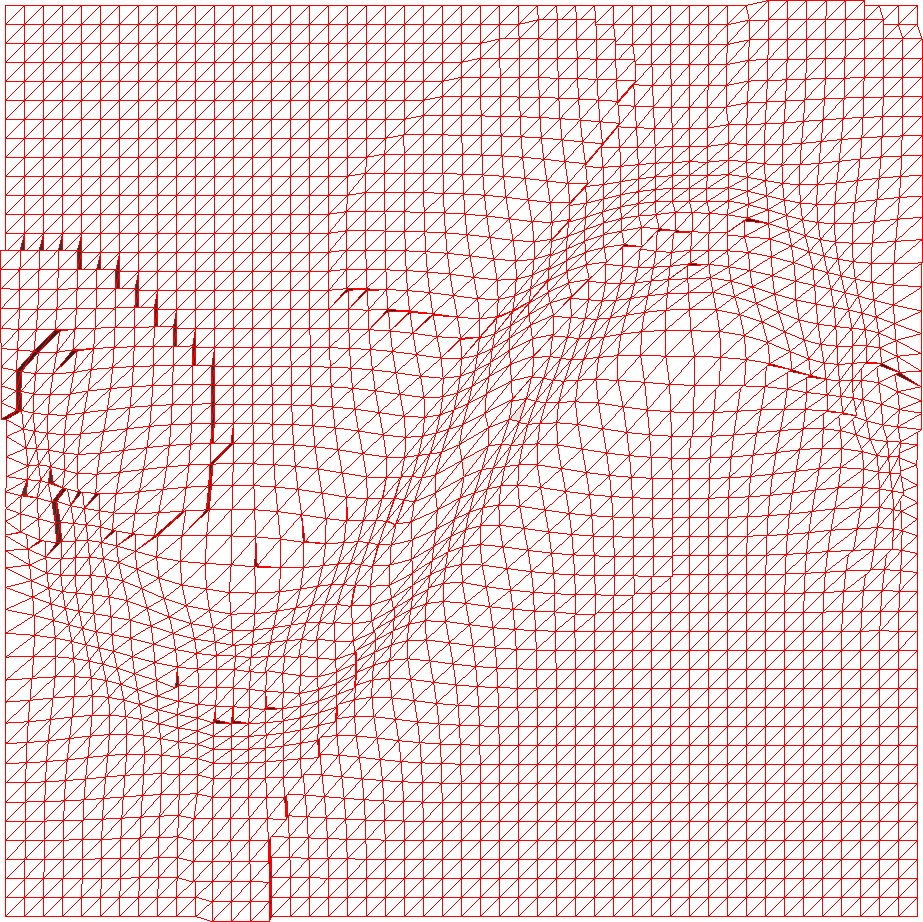}%
      \caption{Winslow's \cref{win-2}}\label{fig:ex:4:winslow}
   \end{subfigure}%
   \hfill%
   \begin{subfigure}[t]{0.31\linewidth}
      \includegraphics[clip, width=1.0\linewidth]{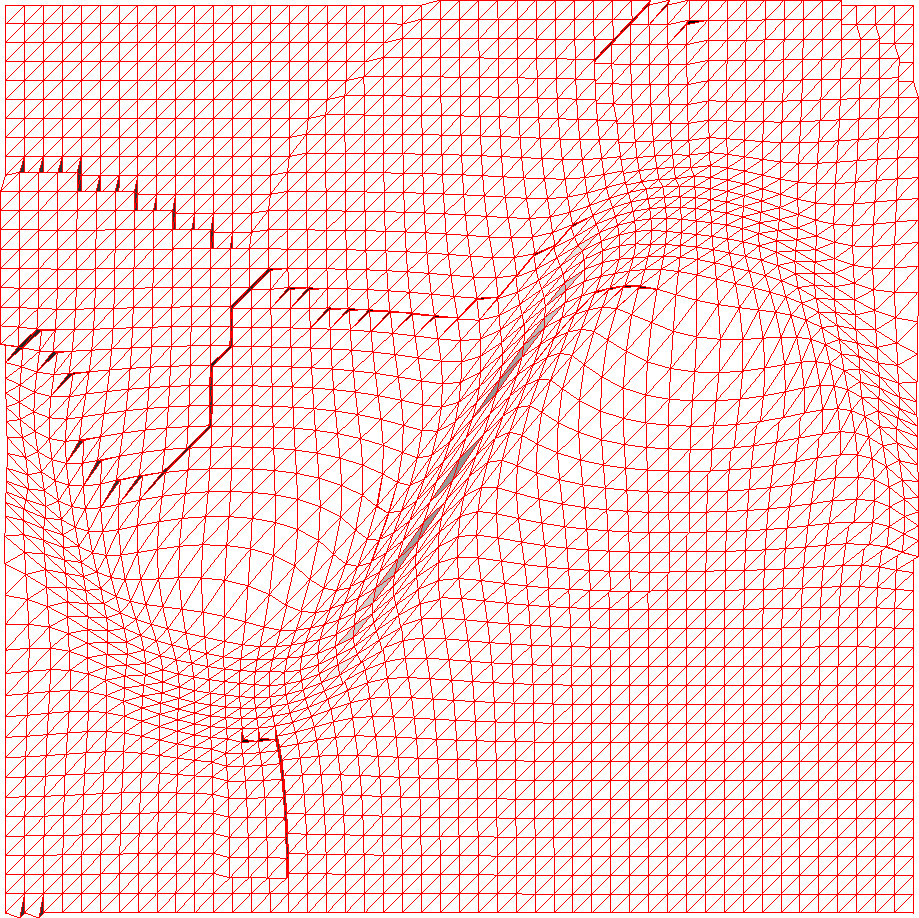}%
      \caption{Huang's \cref{huang-1}}\label{fig:ex:4:huang}
   \end{subfigure}%
   \hfill%
   \begin{subfigure}[t]{0.31\linewidth}
      \includegraphics[clip, width=1.0\linewidth]{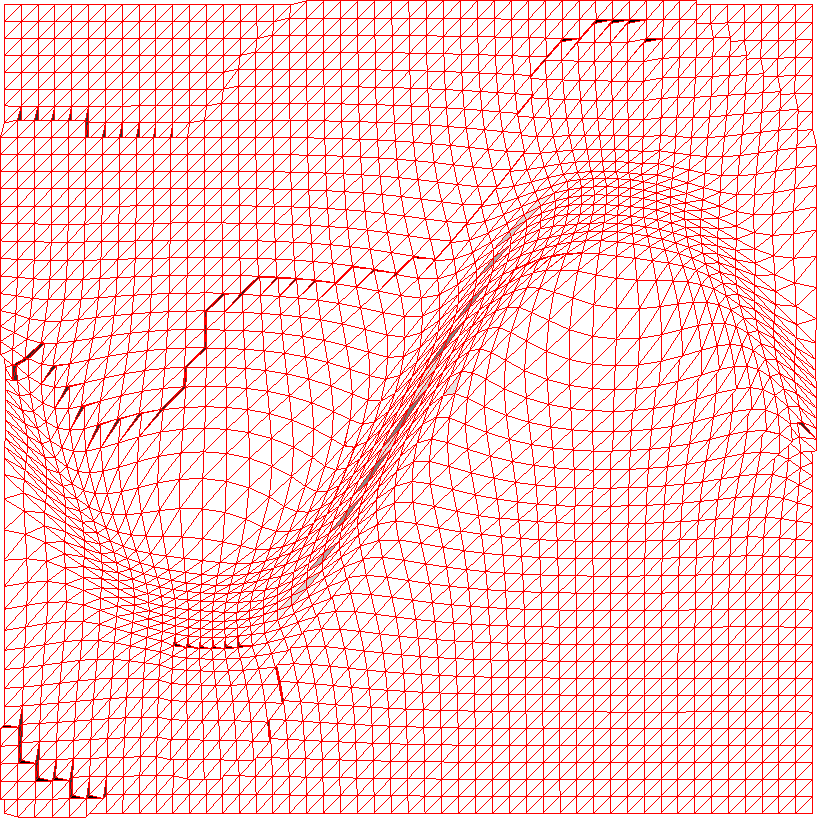}%
      \caption{Huang and Russell's \cref{hr-1}}\label{fig:ex:4:HR}
   \end{subfigure}%
   \\[2.0ex]
   \begin{subfigure}[t]{0.31\linewidth}
      \includegraphics[clip, width=1.0\linewidth]{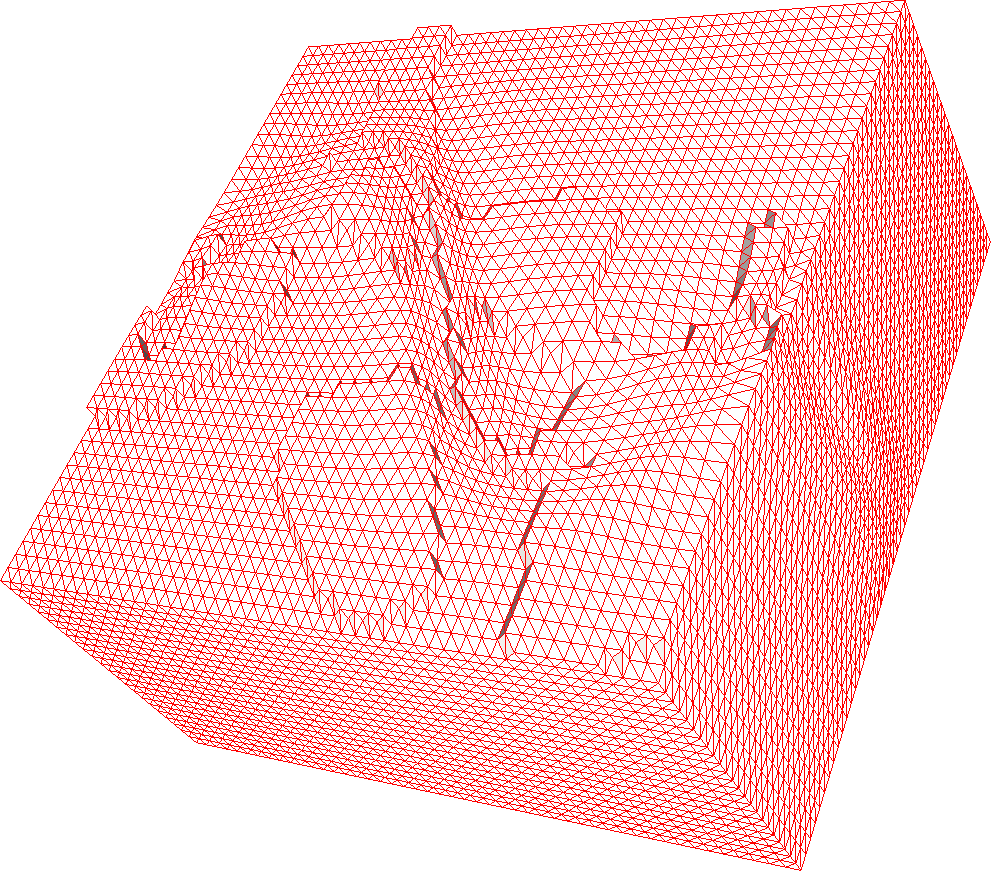}%
      \caption{Winslow's \cref{win-2}}\label{fig:ex:4:winslow+3d}
   \end{subfigure}%
   \hfill%
   \begin{subfigure}[t]{0.31\linewidth}
      \includegraphics[clip, width=1.0\linewidth]{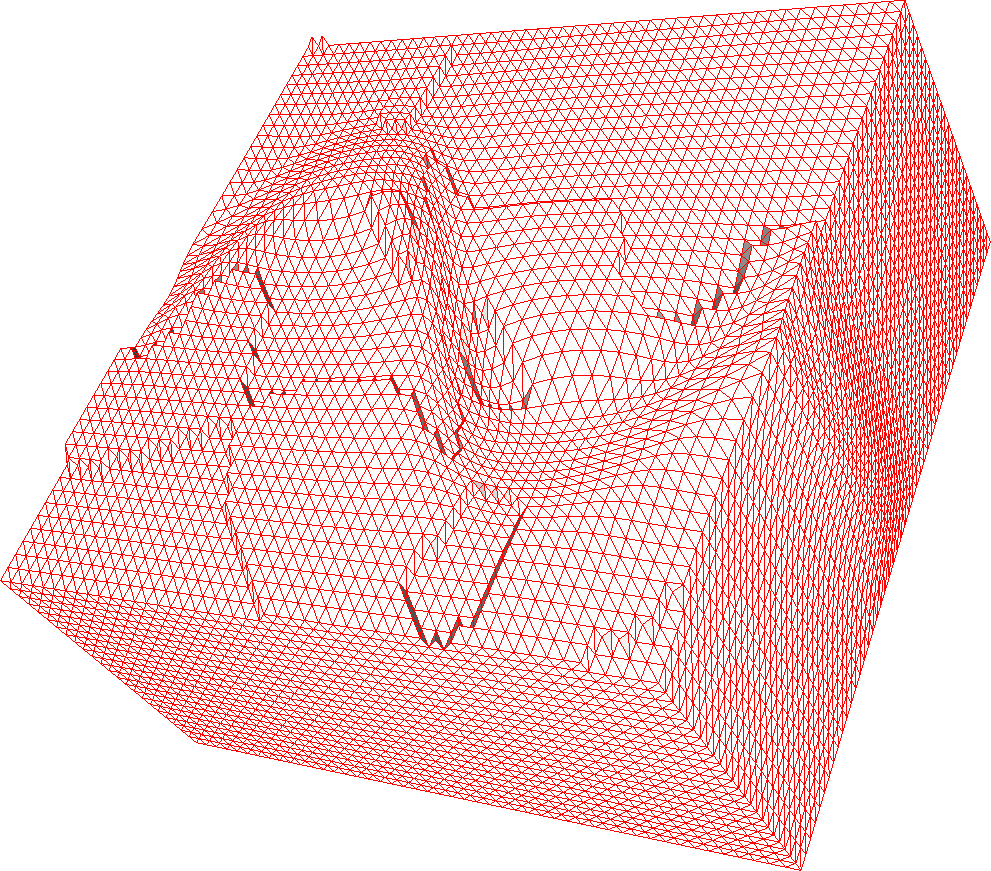}%
      \caption{Huang's \cref{huang-1}}\label{fig:ex:4:huang+3d}
   \end{subfigure}%
   \hfill%
   \begin{subfigure}[t]{0.31\linewidth}
      \includegraphics[clip, width=1.0\linewidth]{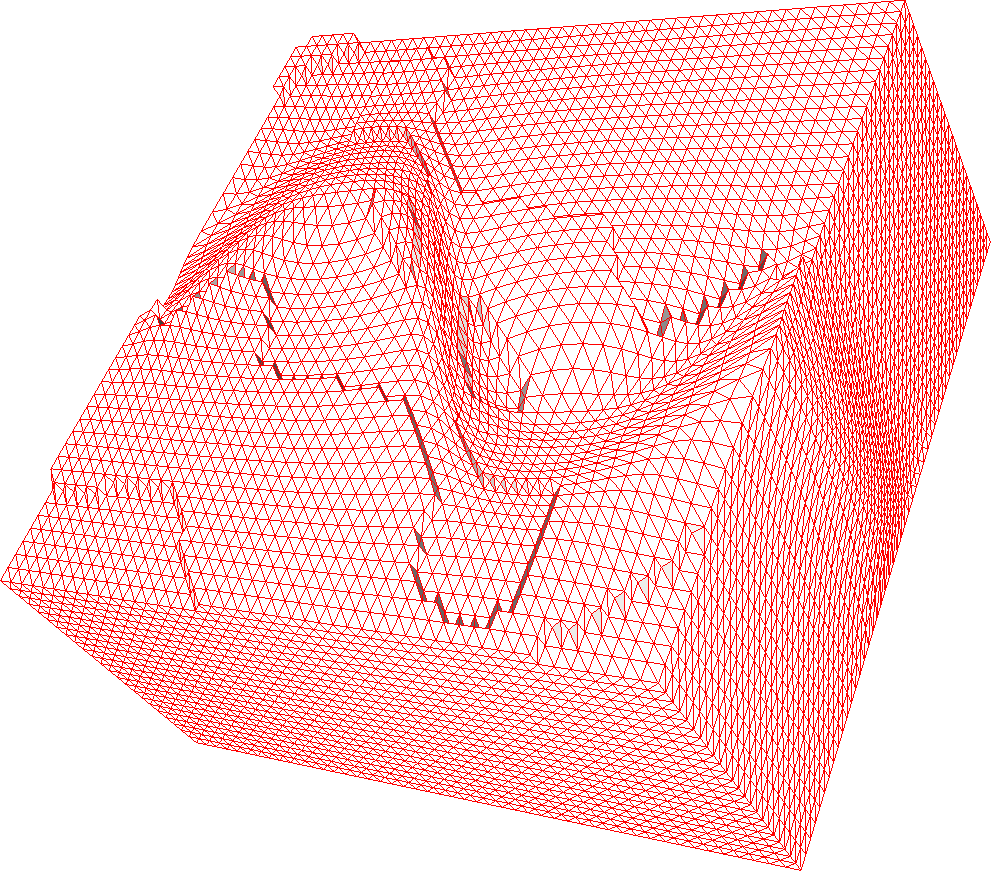}%
      \caption{Huang and Russell's \cref{hr-1}}\label{fig:ex:4:HR+3d}
   \end{subfigure}%
   \\[2.0ex]
   \begin{subfigure}[t]{0.31\linewidth}
      \begin{tikzpicture}
         \begin{semilogxaxis}[
            width=0.75\linewidth,
            height=0.60\linewidth,
            ymin = 0.9, ymax = 2.8, 
         ]
            \addplot[color=\cWinslow, \ltWinslow, mark = +, ]
            table [x index=1, y index=4, col sep = space] {results-4-winslow.dat};
            \addlegendentry{Winslow}
            \addplot[color=\cHuang, \ltHuang, mark = +,]
            table [x index=1, y index=4, col sep = space] {results-4-huang.dat};
            \addlegendentry{Huang}
            \addplot[color=\cHR, \ltHR, mark = +,]
            table [x index=1, y index=4, col sep = space] {results-4-hr.dat};
            \addlegendentry{HR}
         \end{semilogxaxis}
      \end{tikzpicture}
      \caption{$Q_{eq}$ vs. $N$}\label{fig:ex:4:qeq}
   \end{subfigure}%
   \hfill%
   \begin{subfigure}[t]{0.31\linewidth}
        \begin{tikzpicture}
         \begin{semilogxaxis}[
            width=0.75\linewidth,
            height=0.60\linewidth,
            ymin = 0.9, ymax = 2.8, 
         ]
            \addplot[color=\cWinslow, \ltWinslow, mark = +,]
            table [x index=1, y index=5, col sep = space] {results-4-winslow.dat};
            \addlegendentry{Winslow}
            \addplot[color=\cHuang, \ltHuang, mark = +,]
            table [x index=1, y index=5, col sep = space] {results-4-huang.dat};
            \addlegendentry{Huang}
            \addplot[color=\cHR, \ltHR, mark = +,]
            table [x index=1, y index=5, col sep = space] {results-4-hr.dat};
            \addlegendentry{HR}
         \end{semilogxaxis}
      \end{tikzpicture}
      \caption{$Q_{ali}$ vs. $N$}\label{fig:ex:4:qali}
   \end{subfigure}%
   \hfill%
   \begin{subfigure}[t]{0.31\linewidth}
      \begin{tikzpicture}
         \begin{loglogaxis}[
            width=0.75\linewidth,
            height=0.60\linewidth,
            legend style={at={(0.98,1.00)}, anchor=north east,align=left},
         ]
            \addplot[color=\cWinslow, \ltWinslow, mark = +,]
            table [x index=1, y index=2, col sep = space] {results-4-winslow.dat};
            \addlegendentry{Winslow}
            \addplot[color=\cHuang, \ltHuang, mark = +,]
            table [x index=1, y index=2, col sep = space] {results-4-huang.dat};
            \addlegendentry{Huang}
            \addplot[color=\cHR, \ltHR, mark = +,]
            table [x index=1, y index=2, col sep = space] {results-4-hr.dat};
            \addlegendentry{HR}
            \addplot[color=\cN, \ltN, mark = none] coordinates
            {(1000, 0.109) (4000, 0.043) (16000, 0.017) (100000, 0.005)};
        \end{loglogaxis}
      \end{tikzpicture}
      \caption{$L^2$ error vs. $N$}\label{fig:ex:4:error}
   \end{subfigure}%
   \caption{\Cref{ex:4}. The top row: slice cuts  of the meshes.
   The middle row: clip cuts of the meshes. The bottom row: $Q_{eq}$, $Q_{ali}$, and the $L^2$ norm of the linear interpolation error  (black line represents $N^{-2/3}$).}\label{fig:ex:4}%
\end{figure}

\begin{figure}[p]
   \begin{subfigure}[t]{0.31\linewidth}
      \includegraphics[clip, width=1.0\linewidth]{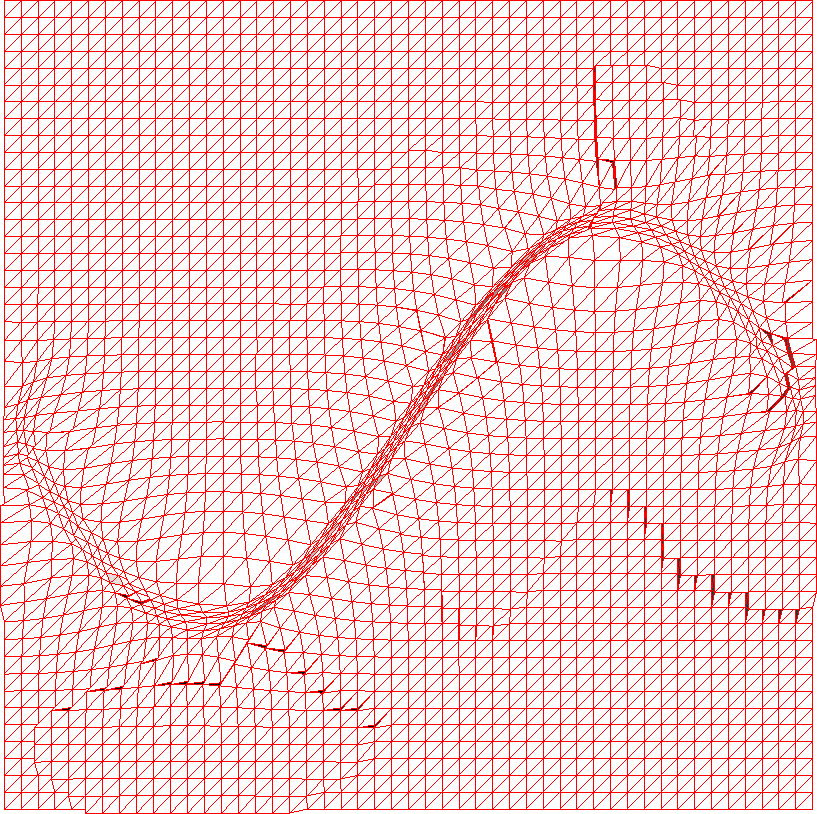}%
      \caption{Winslow's \cref{win-2}}\label{fig:ex:5:winslow}
   \end{subfigure}%
   \hfill%
   \begin{subfigure}[t]{0.31\linewidth}
      \includegraphics[clip, width=1.0\linewidth]{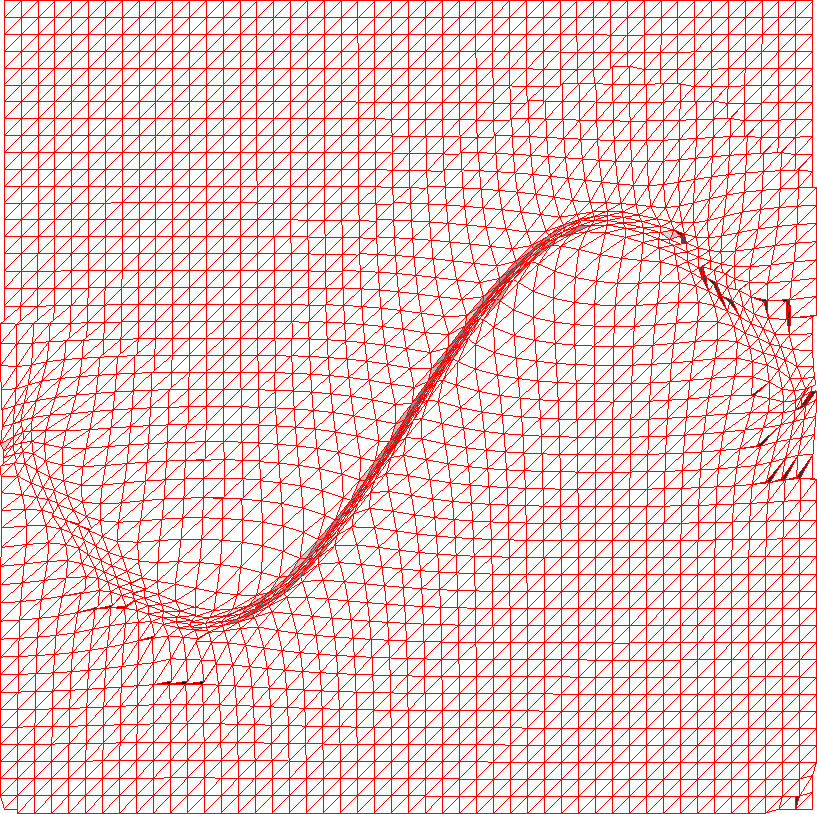}%
      \caption{Huang's \cref{huang-1}}\label{fig:ex:5:huang}
   \end{subfigure}%
   \hfill%
   \begin{subfigure}[t]{0.31\linewidth}
      \includegraphics[clip, width=1.0\linewidth]{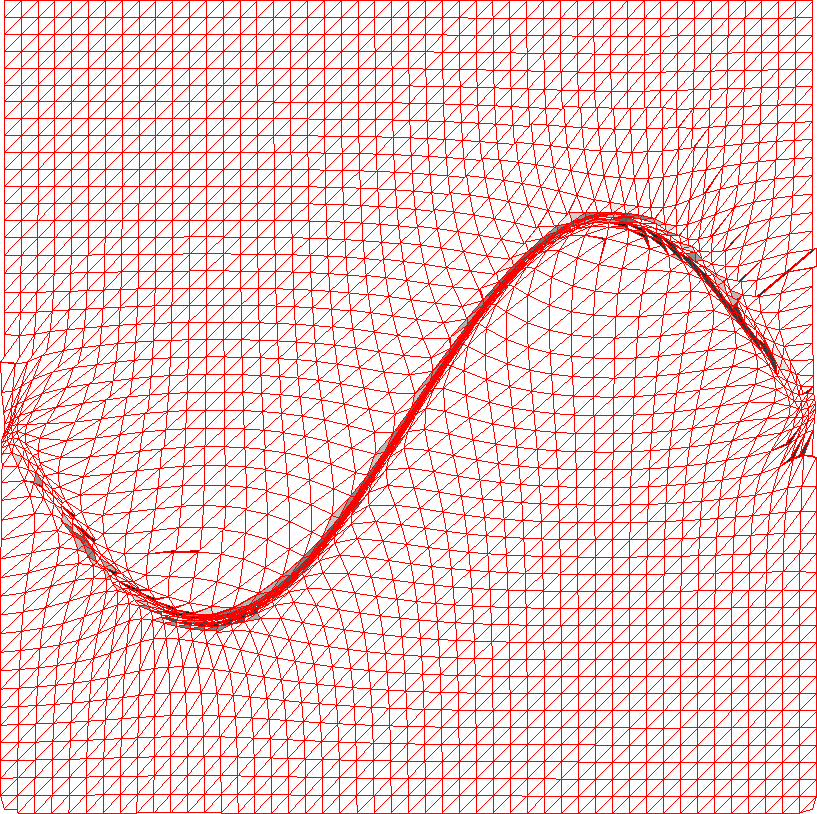}%
      \caption{Huang and Russell's \cref{hr-1}}\label{fig:ex:5:HR}
   \end{subfigure}%
   \\[2.0ex]
   \begin{subfigure}[t]{0.31\linewidth}
      \includegraphics[clip, width=1.0\linewidth]{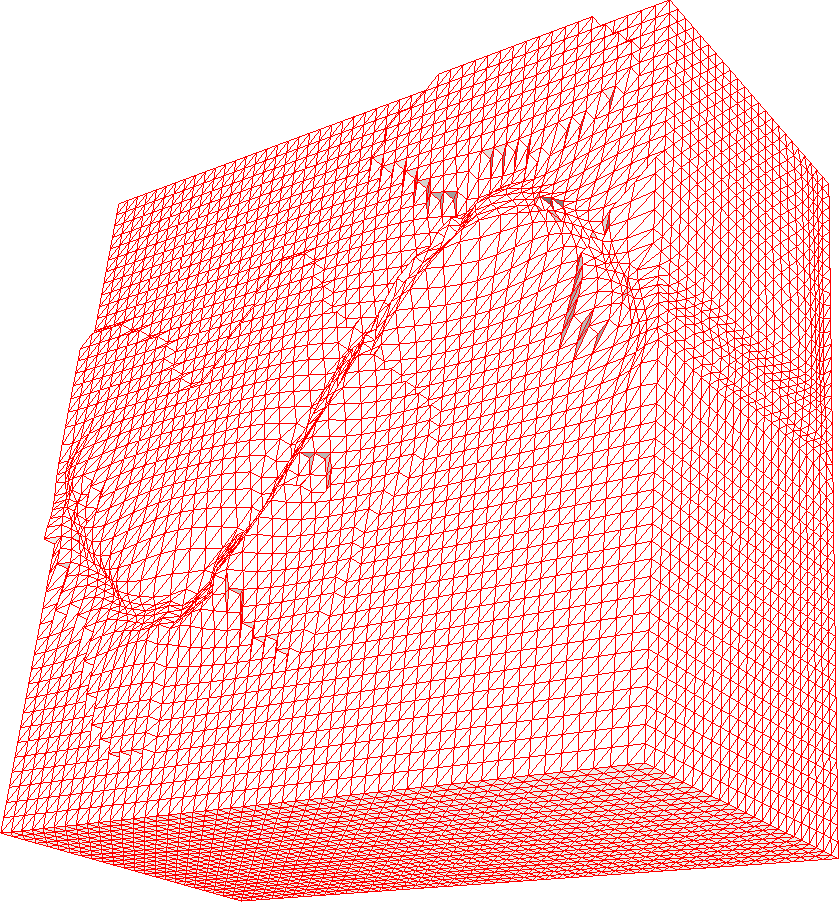}%
      \caption{Winslow's \cref{win-2}}\label{fig:ex:5:winslow+3d}
   \end{subfigure}%
   \hfill%
   \begin{subfigure}[t]{0.31\linewidth}
      \includegraphics[clip, width=1.0\linewidth]{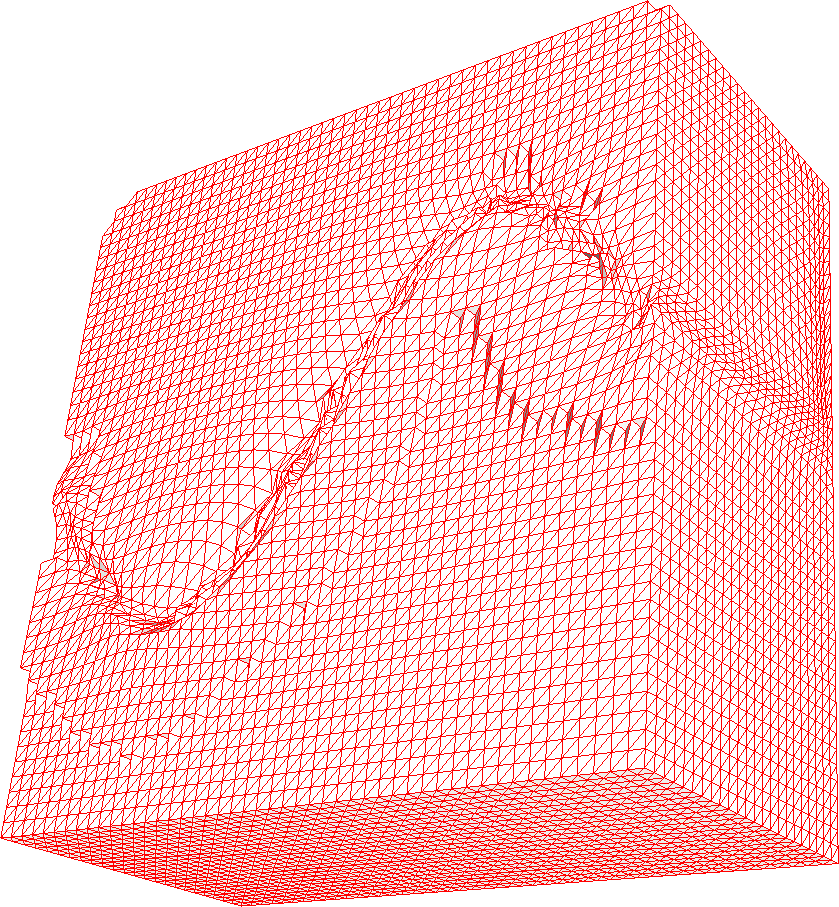}%
      \caption{Huang's \cref{huang-1}}\label{fig:ex:5:huang+3d}
   \end{subfigure}%
   \hfill%
   \begin{subfigure}[t]{0.31\linewidth}
      \includegraphics[clip, width=1.0\linewidth]{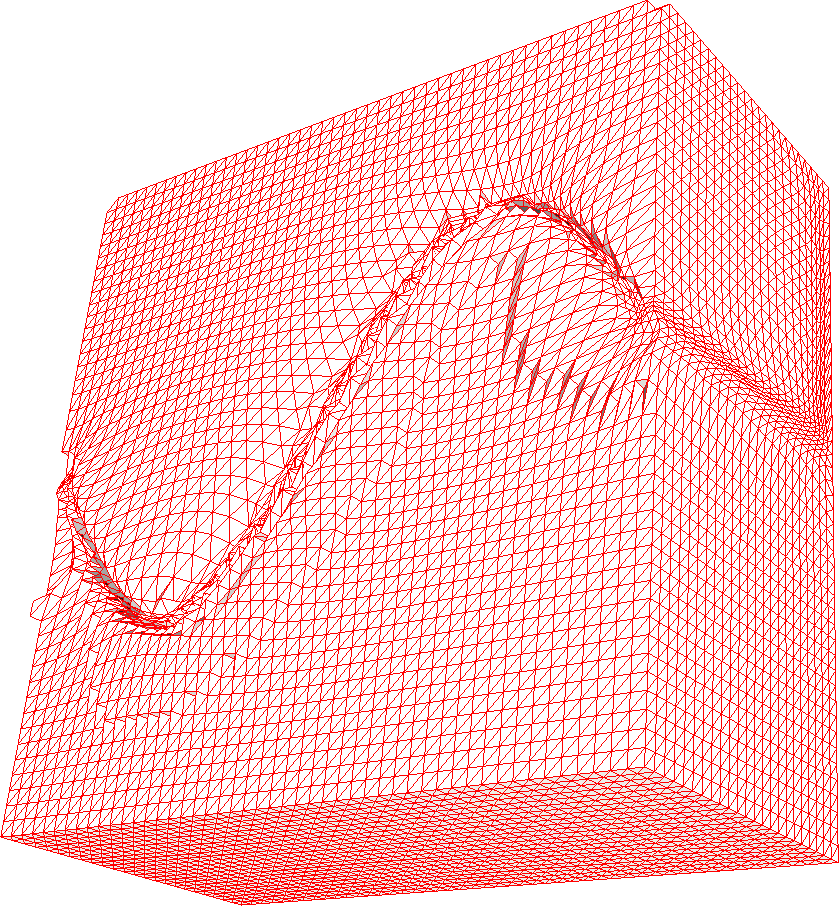}%
      \caption{Huang and Russell's \cref{hr-1}}\label{fig:ex:5:HR+3d}
   \end{subfigure}%
   \\[2.0ex]
   \begin{subfigure}[t]{0.31\linewidth}
      \begin{tikzpicture}
         \begin{semilogxaxis}[
            width=0.75\linewidth,
            height=0.60\linewidth,
            ymin = 0.9, ymax = 2.8, 
         ]
            \addplot[color=\cWinslow, \ltWinslow, mark = +,]
            table [x index=1, y index=4, col sep = space] {results-5-winslow.dat};
            \addlegendentry{Winslow}
            \addplot[color=\cHuang, \ltHuang, mark = +,]
            table [x index=1, y index=4, col sep = space] {results-5-huang.dat};
            \addlegendentry{Huang}
            \addplot[color=\cHR, \ltHR, mark = +,]
            table [x index=1, y index=4, col sep = space] {results-5-hr.dat};
            \addlegendentry{HR}
         \end{semilogxaxis}
      \end{tikzpicture}
      \caption{$Q_{eq}$ vs. $N$}\label{fig:ex:5:qeq}
   \end{subfigure}%
   \hfill%
   \begin{subfigure}[t]{0.31\linewidth}
        \begin{tikzpicture}
         \begin{semilogxaxis}[
            width=0.75\linewidth,
            height=0.60\linewidth,
            ymin = 0.9, ymax = 2.8, 
         ]
            \addplot[color=\cWinslow, \ltWinslow, mark = +,]
            table [x index=1, y index=5, col sep = space] {results-5-winslow.dat};
            \addlegendentry{Winslow}
            \addplot[color=\cHuang, \ltHuang, mark = +,]
            table [x index=1, y index=5, col sep = space] {results-5-huang.dat};
            \addlegendentry{Huang}
            \addplot[color=\cHR, \ltHR, mark = +,]
            table [x index=1, y index=5, col sep = space] {results-5-hr.dat};
            \addlegendentry{HR}
         \end{semilogxaxis}
      \end{tikzpicture}
      \caption{$Q_{ali}$ vs. $N$}\label{fig:ex:5:qali}
   \end{subfigure}%
   \hfill%
   \begin{subfigure}[t]{0.31\linewidth}
      \begin{tikzpicture}
         \begin{loglogaxis}[
            width=0.75\linewidth,
            height=0.60\linewidth,
            legend style={at={(0.98,1.00)}, anchor=north east,align=left},
         ]
            \addplot[color=\cWinslow, \ltWinslow, mark = +,]
            table [x index=1, y index=2, col sep = space] {results-5-winslow.dat};
            \addlegendentry{Winslow}
            \addplot[color=\cHuang, \ltHuang, mark = +,]
            table [x index=1, y index=2, col sep = space] {results-5-huang.dat};
            \addlegendentry{Huang}
            \addplot[color=\cHR, \ltHR, mark = +,]
            table [x index=1, y index=2, col sep = space] {results-5-hr.dat};
            \addlegendentry{HR}
            \addplot[color=\cN, \ltN, mark = none] coordinates
            {(2000, 0.218) (8000, 0.086) (32000, 0.034)};
        \end{loglogaxis}
      \end{tikzpicture}
      \caption{$L^2$ error vs.\ $N$}\label{fig:ex:5:error}
   \end{subfigure}%
   \caption{\Cref{ex:5}. The top row: slice cuts  of the meshes.
   The middle row: clip cuts of the meshes. The bottom row: $Q_{eq}$, $Q_{ali}$, and the $L^2$ norm of the linear interpolation error  (black line represents $N^{-2/3}$).}\label{fig:ex:5}%
\end{figure}
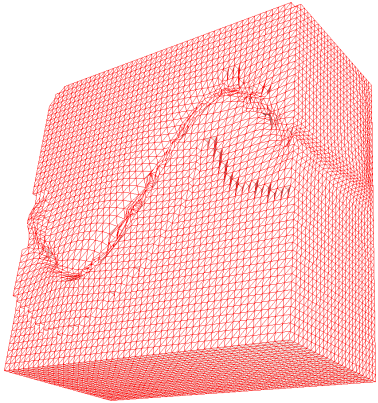
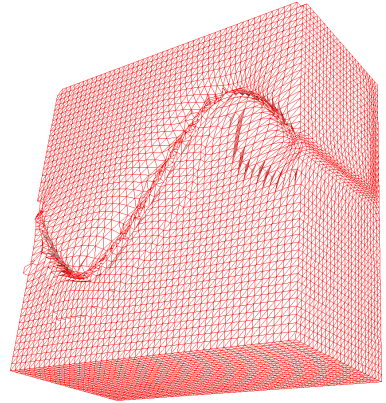
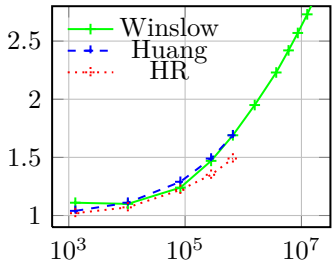
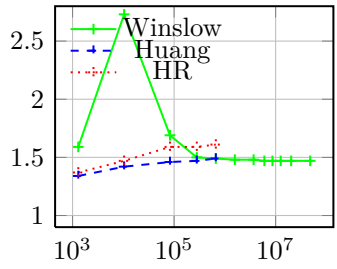
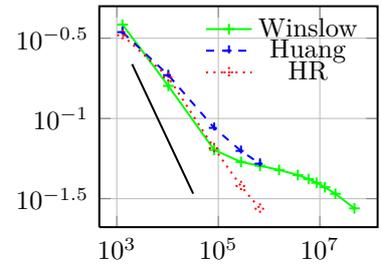

Adaptive mesh examples (slice and clip cuts) and numerical results are given in \cref{fig:ex:3,fig:ex:4,fig:ex:5}.

As in two dimensions, all three functionals provide good size and shape adaptation, with $Q_{eq}$ and $Q_{ali}$ being reasonably small.
The best mesh size control $Q_{eq}$ is given by HR functional (Figs.~\labelcref{fig:ex:3:qeq}, \labelcref{fig:ex:4:qeq}, and \labelcref{fig:ex:5:qeq}), although for the considered examples, HR has a slightly worse mesh alignment quality $Q_{ali}$ than the others (Figs.~\labelcref{fig:ex:3:qali}, \labelcref{fig:ex:4:qali}, and \labelcref{fig:ex:5:qali}).

A closer look at the example meshes (slice cuts) reveals that, as in 2D, W functional ---based on variable diffusion--- is noticeably more aggressive in moving nodes toward the steep features or, alternatively, one can say that the functionals \cref{huang-1,hr-1} based on equidistribution and alignment distribute the nodes with the better correspondence with the given $\M$.
For coarse meshes, all three functionals provide similar results (see convergence plots in Figs.~\labelcref{fig:ex:3:error}, \labelcref{fig:ex:4:error}, \labelcref{fig:ex:5:error}); however, for fine meshes, sizing of mesh elements obtained by means of W functional is not quite as good as for H and HR functionals, as indicated by a larger $Q_{eq}$.

Altogether, the linear interpolation error (Figs.~\labelcref{fig:ex:3:error}, \labelcref{fig:ex:4:error}, and \labelcref{fig:ex:5:error}) suggests that HR functional provides the best mesh, followed by H and W functionals.
One may notice from \cref{fig:ex:4:error,fig:ex:5:error} that the convergence of the linear interpolation error for W functional slows down near $N= 10^5$ for \cref{ex:4,ex:5}, although it seems to improve as the mesh is refined (\cref{fig:ex:5:error}).
The reason for this behaviour is not clear to us.
On the other hand, W functional has the simplest form and seems to be more economic to compute than the other two.
From tentative comparison, mesh generation using W functional uses about one fifth to an half of the CPU time used with H or HR functional.
Qualitatively,  this is not difficult to understand since W functional is convex whereas the others are not (although they are polyconvex).

\section{Conclusions}
\label{SEC:conclusions}

Among the three functionals in this study, Huang and Russell's functional consistently provides the best mesh for piecewise linear interpolation in both two and three dimensions.
In all examples it leads to the best equidistribution quality and the smallest interpolation error.
Interestingly, while it results in the best mesh alignment quality in two dimensions, the functional gives a slightly worse mesh alignment than the other two functionals in three dimensions. 

Meshes obtained by means of Winslow's functional have the worst mesh equidistribution (element size) quality and the largest interpolation error in four out of five examples, although in three dimensions its mesh alignment is quite good and even better than that of the meshes obtained using Huang and Russell's functional.
An explanation to this behavior could be the fact that this functional does not have an explicit mechanism to control the equidistribution property.

The behavior of Huang's functional is somewhere in between Winslow's and Huang and Russell's functionals: both in mesh quality measures and interpolation error.
It provides better mesh sizing than Winslow's functional but not quite as good as Huang and Russell's.
On the other hand, it provides the best (or very close to the best) mesh alignment in all examples.

While being able to produce correct and good quality mesh concentration, Winslow's functional seems to have the tendency to move more points toward the area of interest and is slightly less reliable than the other two functionals especially when the mesh is fine.
On the other hand, it has a very simple form and is more economic to compute than the others.
It can be a good choice for mesh adaptation at least for coarser meshes, for which all of the three functionals produce comparable meshes.

Finally, it should be pointed out that the numerical experiment we conducted in this work is limited and more work is needed to have an extensive and more complete understanding of the behavior of the meshing functionals especially in three dimensions.
Moreover, the newly developed implementation of the variational methods in~\cite{HK2014} has been crucial to the current study to perform substantial computations in two and three dimensions.
 It is our hope that it can serve as an efficient tool for use in future studies of mesh adaptation and movement.


\end{document}